\documentclass[11pt,twoside,epsf]{article}
\usepackage{amsmath,latexsym,amssymb,amsfonts,amsbsy}
\usepackage{bm}
\usepackage{graphicx,subfigure}
\usepackage{cases}
\newfont{\bb}{msbm10}
\def\Bbb#1{\mbox{\bb #1}}

\def\opt{{\rm opt}}

\def\zspace{{\rm null}}

\def\sp{{\rm sp}}

\newtheorem{method}{Method}[section]

\newtheorem{remark}{Remark}[section]
\newtheorem{theorem}{Theorem}[section]

\newtheorem{lemma}{Lemma}[section]
\newtheorem{corollary}{Corollary}[section]
\newtheorem{assumption}{Assumption}[section]

\newcommand{\ceals}{\makebox{{\Bbb C}}}

\newcommand{\bbb}[1]{\text{\bf #1}}
\newcommand{\txt}[1]{\text{\rm #1}}
\newcommand{\seq}[2]{#1_{\text{\tiny$\triangle t$}}^{#2}}
\newcommand{\seqi}[3]{#1_{\text{\tiny$\triangle t,{#2}$}}^{#3}}

\makeatletter 

\@addtoreset{equation}{section}
\makeatother 

\begin{document}
\cleardoublepage
\pagestyle{myheadings}

\bibliographystyle{plain}

\title{Discrete-Time\\ Accelerated Block Successive
       Overrelaxation Methods\\
       for Time-Dependent Stokes Equations
\footnote{Supported by The National Natural Science Foundation
(No. 11101213), P.R. China.}
}
\author
{Xi Yang
\footnote{Corresponding author at: Department of Mathematics,
Nanjing University of Aeronautics and Astronautics,
No. 29 Yudao Street, Nanjing 210007, P.R. China.
Email: yangxi@lsec.cc.ac.cn}\\
\\
{\it State Key Laboratory of Scientific/Engineering Computing}\\
{\it Institute of Computational Mathematics and
        Scientific/Engineering Computing}\\
{\it Academy of Mathematics and Systems Science}\\
{\it Chinese Academy of Sciences, P.O. Box 2719,
Beijing 100190, P.R. China}\\
\\
{\it Department of Mathematics}\\
{\it Nanjing University of Aeronautics and Astronautics}\\
{\it No. 29 Yudao Street, Nanjing 210007, P.R. China}\\
}

\maketitle

\markboth{\small X. Yang}
{\small DABSOR Methods for Time-Dependent Stokes Equations}

\begin{abstract}
To further study the application of waveform relaxation methods in fluid dynamics
in actual computation, this paper provides a general theoretical analysis of
discrete-time waveform relaxation methods for solving linear DAEs. A class of
discrete-time waveform relaxation methods, named discrete-time accelerated block
successive overrelaxation (DABSOR) methods, is proposed for solving linear DAEs
derived from discretizing time-dependent Stokes equations in space by using
``Method of Lines''. The analysis of convergence property and optimality of
the DABSOR method are presented in detail. The theoretical results and the
efficiency of the DABSOR method are verified by numerical experiments.

\bigskip


\bigskip

{\bf AMS(MOS) Subject Classifications: } 65F10, 65L80,
65N06, 65N22, 65N40; CR: G1.3.

\end{abstract}

\section{Introduction}
\label{sec-introduction}
We consider the numerical solution for time-dependent Stokes
equations of the form
\begin{eqnarray}
\label{time-stokes}
\left\{
    \begin{array}{rcll}
    \frac{\partial \vec {u}}{\partial t}
    -\nu\nabla^2\vec {u} + \nabla p &=&
    0 &\mbox{in}\enskip\Omega\times\mathbb{R}_+,\\
    \nabla\cdot\vec {u}&=&
    0 &\mbox{in}\enskip\Omega\times\mathbb{R}_+,\\
    \vec {u}&=&0 &
    \mbox{on}\enskip\partial\Omega\times\mathbb{R}_+,\\
    \vec {u}&=&\vec {u}_0 &
    \mbox{on}\enskip\Omega\times\{0\},\\
    \end{array}
    \right.
\end{eqnarray}
with certain initial and boundary conditions in $d$-dimensional
locally Lipschitz bounded domain $\Omega\subset\mathbb{R}^d$
($d=2$ or $3$). By adopting the idea of ``Method of Lines'',
equations (\ref{time-stokes}) are discretized in space to obtain
the following saddle-structured differential-algebraic equations
(DAEs)
\begin{eqnarray}
\label{dis-dae-LCC}
\mathcal{B}\dot z(t) + \mathcal{A} z(t) = b(t),
\quad z(0) = z_{0},
\end{eqnarray}
where $z(t)=\left(x(t)^T, y(t)^T\right)^T$ with
$x(t)$ and $y(t)$ related to velocity and pressure in
equations (\ref{time-stokes}) respectively,
$\mathcal{B}$ and $\mathcal{A}$ are block two-by-two square
matrices of the form
$$
\mathcal{B}=\left(\begin{array}{cc}
I&0\\0&0
\end{array}\right)\quad \mbox{and} \quad
\mathcal{A}=\left(\begin{array}{cc}
A&B\\-B^{T}&0
\end{array}\right),
$$
with $A\in\mathbb{R}^{r\times r}$ being a symmetric
positive definite matrix,
$B\in\mathbb{R}^{r\times l}$ a full column-rank
matrix, and $I\in\mathbb{R}^{r\times r}$
the identity matrix. Here, $r$ and $l$ are known
positive integers.

In \cite{BaiYang12}, a class of continuous-time waveform relaxation
methods for solving linear DAEs (\ref{dis-dae-LCC}) has been proposed
by the application of generalized successive overrelaxation (GSOR)
technique \cite{gsor}. Previous continuous-time waveform relaxation methods,
for solving ODEs
\cite{wr-lcc-4,JAN96,wr-sor,wr-lcc-1,Nev89II,wr-lcc-3,wr-lcc-5,WangBai05}
and DAEs \cite{BaiYang12,BaiYang11,JiangW-00,wr-lcc-2}, can be regarded as
extensions of the classical iterative methods for solving system of
algebraic equations with iterating space changing from $\mathbb{R}^n$ to
the waveform space. Since the analytical operations and exact expressions
of waveforms are required in each iterative step, the continuous-time
waveform relaxation methods are unlikely to be practical
numerical methods, but only be of theoretical interest.

In actual numerical solution for linear DAEs (\ref{dis-dae-LCC}), the
continuous-time waveform relaxation methods in \cite{BaiYang12} are replaced by
a class of discrete-time waveform relaxation methods,
named {\em discrete-time accelerated block successive overrelaxation}
({\bf DABSOR}) methods. This paper studies the general
theory of discrete-time waveform relaxation methods for solving linear DAEs
and the special case of the
DABSOR method. In previous discrete-time waveform relaxation
methods \cite{JAN96,Jiang04,Jiang09}, the waveforms and the linear differential
operators in continuous-time waveform relaxation methods are replaced
by vector sequences and discrete linear convolution operators respectively.
Here, the vector sequences are composed of values of waveforms on each time
level, and the discrete linear convolution operators are closely related to
the discretizations of the linear differential operators by different
time stepping schemes \cite{LR-Petzold, G-Wanner} in each iterative
step of the continuous-time waveform relaxation methods for solving
linear DAEs (\ref{dis-dae-LCC}).

The structure of this paper is as follows. It is started in Section \ref{sec-dis-conv}
by briefly reviewing the spectral properties of the discrete linear convolution operator.
The general framework of the discrete-time waveform relaxation methods and the corresponding
discrete linear convolution operator form are stated in Section \ref{sec-dis-framework}.
In Section \ref{sec-dis-con-ana}, the convergence properties of the discrete linear
convolution operator derived from the discrete-time waveform relaxation methods for solving
linear DAEs are
analyzed on both finite and infinite time interval. The DABSOR method without or with
windowing technique is constructed
in Section \ref{sec-dabsor}, and the convergence domain of relaxation parameters and
the optimality of the DABSOR method are also presented here. The numerical results
are listed in Section \ref{sec-num}, which is followed by the concluding remarks in
Section \ref{sec-cremarks}.

\section{Spectral Properties of Discrete Linear Convolution Operator}
\label{sec-dis-conv}
Consider the following iterative scheme
\begin{eqnarray}
\label{dis-conv-op}
z_{\text{\tiny$\triangle t$}}^{(k)}=
\mathcal{H}_{\text{\tiny$\triangle t$}}
z_{\text{\tiny$\triangle t$}}^{(k-1)}
+\varphi_{\text{\tiny$\triangle t$}},
\end{eqnarray}
where the subscript $\triangle t$ is the notation of vector
sequences, e.g. $z_{\text{\tiny$\triangle t$}}^{(k)} =
\{z_i^{(k)}\}_{i=0}^{L-1}$, where $L$ (possibly infinite) is the
number of components. Each $d$-dimensional component denotes
the approximate solution of a $d$-dimensional ODEs or DAEs on
a time level. Operator
$\mathcal{H}_{\text{\tiny$\triangle t$}}$ is a discrete linear
convolution operator with matrix-valued kernel
$h_{\text{\tiny$\triangle t$}}$,
\begin{eqnarray}
(\mathcal{H}_{\text{\tiny$\triangle t$}}
z_{\text{\tiny$\triangle t$}})_j=
(h_{\text{\tiny$\triangle t$}}\star
z_{\text{\tiny$\triangle t$}})_j=
\sum_{i=0}^j h_{j-i}z_i,\quad j=0,\ldots,L-1.
\end{eqnarray}
The convergence properties of operator
$\mathcal{H}_{\text{\tiny$\triangle t$}}$ are analyzed
in the Banach spaces of $\mathbb{C}^d$-valued $p$-summable sequences
of length $L$, $l_p(L;\mathbb{C}^d)$, or $l_p(L)$ for short,
with norms given by
\begin{eqnarray}
\|z_{\text{\tiny$\triangle t$}}\|_p = \left\{
\begin{array}{ll}
    \root p \of {\sum _{i=0}^{L-1}\|z_i\|^p} & 1\le p<\infty,\\
    \sup \limits_{0\le i <L}\{\|z_i\|\} & p = \infty,
\end{array}
\right.
\end{eqnarray}
with $\|\cdot\|$ any prescribed $\mathbb{C}^d$ vector norm. It is known
that the iterative scheme (\ref{dis-conv-op}) is convergent if and only
if the spectral radius of the discrete linear convolution operator
$\mathcal{H}_{\text{\tiny$\triangle t$}}$, denoted by
$\rho(\mathcal{H}_{\text{\tiny$\triangle t$}})$, is smaller than one.
The following two lemmas provide specific descriptions of
$\rho(\mathcal{H}_{\text{\tiny$\triangle t$}})$ on both finite
and infinite time intervals \cite{JAN96}.
\begin{lemma}
\label{dis-spe-fin}
Consider $\mathcal{H}_{\text{\tiny$\triangle t$}}$ as an operator
in $l_p(L)$, with $1\le p\le \infty$ and $L$ finite. Then,
$\mathcal{H}_{\text{\tiny$\triangle t$}}$ is a bounded operator
and
\begin{eqnarray*}
\rho(\mathcal{H}_{\text{\tiny$\triangle t$}})=
\rho(h_0)=\rho(\bbb{H}_{\text{\tiny$\triangle t$}}(\infty)),
\end{eqnarray*}
with $\bbb{H}_{\text{\tiny$\triangle t$}}(s)=
\sum _{i=0}^{L-1}h_is^{-i}$ the discrete Laplace transform of
$h_{\text{\tiny$\triangle t$}}$.
\end{lemma}

\begin{lemma}
\label{dis-spe-inf} Suppose that $h_{\text{\tiny$\triangle
t$}}\in\l_1(\infty)$, and consider
$\mathcal{H}_{\text{\tiny$\triangle t$}}$ as an operator in
$l_p(\infty)$, with $1\le p\le \infty$. Then,
$\mathcal{H}_{\text{\tiny$\triangle t$}}$ is a bounded operator and
\begin{eqnarray*}
\rho(\mathcal{H}_{\text{\tiny$\triangle t$}})=
\max_{|s|\ge 1}\rho(\bbb{H}_{\text{\tiny$\triangle t$}}(s))
=\max_{|s| = 1}\rho(\bbb{H}_{\text{\tiny$\triangle t$}}(s)),
\end{eqnarray*}
with $\bbb{H}_{\text{\tiny$\triangle t$}}(s)=
\sum _{i=0}^{\infty}h_is^{-i}$ the discrete Laplace transform of
$h_{\text{\tiny$\triangle t$}}$.
\end{lemma}

\section{Discrete-Time Waveform Relaxation Methods}
\label{sec-dis-framework}

The continuous-time waveform relaxation methods for
solving the linear DAEs (\ref{dis-dae-LCC}) are defined by
spliting the square matrices $\mathcal{B}$ and
$\mathcal{A}\in \mathbb{R}^{(r+l)\times (r+l)}$ into
\begin{eqnarray*}
\mathcal{B} = M_{\mathcal{B}} - N_{\mathcal{B}}
\quad \mbox{and} \quad
\mathcal{A} = M_{\mathcal{A}} - N_{\mathcal{A}},
\end{eqnarray*}
respectively. Then the corresponding
iterative scheme is of the form
\begin{eqnarray}
\label{wr-fw}
\left\{
    \begin{array}{l}
    M_{\mathcal{B}}\dot z^{(k)}
    +M_{\mathcal{A}}z^{(k)}=N_{\mathcal{B}}\dot
    z^{(k-1)}+N_{\mathcal{A}}z^{(k-1)}+b,\\
    z^{(k)}(0)=z_{0}.
    \end{array}
\right.
\end{eqnarray}
Furthermore, iterative scheme (\ref{wr-fw}) can be rewritten
explicitly
\begin{align}
& z^{(k)} = \mathcal{K}(z^{(k-1)}) + \Phi(b),\\
\intertext{where}
& \mathcal{K}(z) = (\mathcal{L}^{-1}
(sM_{\mathcal{B}}+M_{\mathcal{A}})^{-1}
(sN_{\mathcal{B}}+N_{\mathcal{A}})
\mathcal{L})(z) \nonumber\\
\intertext{and}
& \Phi(b) = (\mathcal{L}^{-1}
(sM_{\mathcal{B}}+M_{\mathcal{A}})^{-1}
\mathcal{L})(b).
\nonumber
\end{align}
Here, $\mathcal{L}$ denotes the continuous Laplace
transform. It has been shown in \cite{BaiYang11} that
\begin{eqnarray}
\rho(\mathcal{K})=\sup _{\Re(s)=\sigma}
\rho\left( \bbb{K}(s) \right),
\end{eqnarray}
where
\begin{eqnarray}
\label{con-lap-mat}
\bbb{K}(s)=(sM_{\mathcal{B}}+M_{\mathcal{A}})^{-1}
(sN_{\mathcal{B}}+N_{\mathcal{A}}).
\end{eqnarray}

Applying a linear multistep formula to the
continuous-time waveform relaxation scheme
(\ref{wr-fw}) leads to the following
discrete-time waveform relaxation scheme
{\setlength{\multlinegap}{50pt}
\begin{multline}
\label{dis-wr-fw}
\frac{1}{\triangle t} \sum _{j=0}^{\nu} \alpha_j
M_{\mathcal{B}} z_{n+j}^{(k)}+\sum _{j=0}^{\nu} \beta_j
M_{\mathcal{A}} z_{n+j}^{(k)}=\\
\frac{1}{\triangle t} \sum _{j=0}^{\nu} \alpha_j
N_{\mathcal{B}} z_{n+j}^{(k-1)}+\sum _{j=0}^{\nu} \beta_j
N_{\mathcal{A}} z_{n+j}^{(k-1)}+
\sum _{j=0}^{\nu}\beta_j b_{n+j},\quad n\ge 0.
\end{multline}}
Assume that the $\nu$ starting values of the linear multistep
formula are known, hence it is not necessary to iterate on the $\nu$ starting
values, i.e. $z_j^{(k)}=z_{j}^{(k-1)}=z_j$, for $j<\nu$.
Due to the numerical stability,
the rest of this paper is concentrated on the application of implicit linear
multistep formulae, i.e. $\beta_{\nu}\ne 0$.

For every nonnegative integer $n$, system of linear equations
(\ref{dis-wr-fw}) can be solved uniquely if and only if the
following condition is satisfied
\begin{eqnarray}
\label{dis-sol-con}
\frac{\alpha_{\nu}}{\beta_{\nu}}\notin
\sp(M_{\mathcal{B}}, -\triangle t M_{\mathcal{A}}),
\end{eqnarray}
where $\sp(\cdot)$ represents the spectrum of the matrix pencil
$(M_{\mathcal{B}}, -\triangle t M_{\mathcal{A}})$.
Subsequently, the above condition (\ref{dis-sol-con}) is
referred to as the {\bf discrete solvability condition}.

Similar to the calculation in \cite{Nev89II}, the iterative
scheme (\ref{dis-wr-fw}) can be rewritten into the following
discrete linear convolution operator form
\begin{eqnarray}
\label{dis-conv-wrop}
z_{\text{\tiny$\triangle t$}}^{(k)}=
\mathcal{K}_{\text{\tiny$\triangle t$}}
z_{\text{\tiny$\triangle t$}}^{(k-1)}
+\varphi_{\text{\tiny$\triangle t$}}.
\end{eqnarray}
Since it does not iterate on the $\nu$ starting values, a slight
change is made on the subscript $\triangle t$ here, that is
\begin{eqnarray}
z_{\text{\tiny$\triangle t$}}^{(k)} =
\{z_{\nu +i}^{(k)}\}_{i=0}^{L-1}.
\end{eqnarray}
In order to analyze the properties of the discrete linear convolution
operator $\mathcal{K}_{\text{\tiny$\triangle t$}}$, the
computational error on the $k$-th iteration of
(\ref{dis-wr-fw}) is denoted by $e_n^{(k)}=z_n^{(k)}-z_n$, where
$z_n$ is the exact solution of the disrete system derived from
the discretezation of linear DAEs (\ref{dis-dae-LCC}) by the
corresponding linear multistep formula. Let
$C_j=\frac{1}{\triangle t}\alpha_j
M_{\mathcal{B}}+\beta_j M_{\mathcal{A}}$ and
$D_j=\frac{1}{\triangle t}\alpha_j
N_{\mathcal{B}}+\beta_j N_{\mathcal{A}}$,
based on (\ref{dis-wr-fw}), we get
\begin{eqnarray}
\label{dis-eqn-err}
\sum _{j=0}^{\nu} C_je_{n+j}^{(k)}=
\sum _{j=0}^{\nu} D_je_{n+j}^{(k-1)},
\quad n\ge 0.
\end{eqnarray}
Combine the first $L$ equations,
introduce vector notation $E^{(k)}=\left({e_{\nu}^{(k)}}^T,
{e_{\nu+1}^{(k)}}^T, \ldots, {e_{L+\nu-1}^{(k)}}^T\right)^T$,
and note that $e_j^{(k)}=e_j^{(k-1)}=0$, $j<\nu$, we have
\begin{eqnarray*}
E^{(k)} = C^{-1}DE^{(k-1)}.
\end{eqnarray*}
Here, matrices $C$ and $D$ are $L\times L$-block lower
triangular Toeplitz matrices with $\nu+1$ constant
diagonals. It follows that matrix $C^{-1}D$ is also
a $L\times L$-block lower triangular Toeplitz matrix
with $\nu+1$ constant diagonals. Hence,
$\mathcal{K}_{\text{\tiny$\triangle t$}}$ is a
discrete linear convolution operator on the Banach space $l_p(L)$.

The discrete Laplace transform of the convolution kernel
$\kappa_{\text{\tiny$\triangle t$}}$ of the
discrete linear convolution operator
$\mathcal{K}_{\text{\tiny$\triangle t$}}$
can be obtained by applying discrete Laplace transform
to equation (\ref{dis-eqn-err}). If
$\tilde e_{\text{\tiny$\triangle t$}}^{(k)}(s)$
denotes the discrete Laplace transform of
$e_{\text{\tiny$\triangle t$}}^{(k)}$,
we find
\begin{eqnarray*}
\tilde e_{\text{\tiny$\triangle t$}}^{(k)}(s)
=\bbb{K}_{\text{\tiny$\triangle t$}}(s)
\tilde e_{\text{\tiny$\triangle t$}}^{(k-1)}(s),
\end{eqnarray*}
with discrete Laplace transform of the convolution kernel
$\kappa_{\text{\tiny$\triangle t$}}$
given by
\begin{eqnarray}
\label{dis-lap-mat}
\bbb{K}_{\text{\tiny$\triangle t$}}(s)=
(a(s)M_{\mathcal{B}}+
\triangle t b(s)M_{\mathcal{A}})^{-1}
(a(s)N_{\mathcal{B}}+
\triangle t b(s)N_{\mathcal{A}}).
\end{eqnarray}
where $a(s)=\sum _{j=0}^{\nu}\alpha_js^j$
and $b(s)=\sum _{j=0}^{\nu}\beta_js^j$.
By comparison to (\ref{con-lap-mat}), we have the following
relation
\begin{eqnarray}
\label{rel-dis-con}
\bbb{K}_{\text{\tiny$\triangle t$}}(s)=
\bbb{K}\left( \frac{1}{\triangle t}
\frac{a}{b}(s) \right).
\end{eqnarray}

\section{Convergence Analysis of $\mathcal{K}_{\text{\tiny$\triangle t$}}$}
\label{sec-dis-con-ana}

The convergence property of the discrete linear
convolution operator $\mathcal{K}_{\text{\tiny$\triangle t$}}$
on finite time interval is an immediate result of Lemma \ref{dis-spe-fin}.
The result can be stated as the following theorem straightforwardly.
\begin{theorem}
\label{th-dis-con-fin}
Assume that the discrete solvability condition (\ref{dis-sol-con}) is satisfied,
and consider $\mathcal{K}_{\text{\tiny$\triangle t$}}$
as a discrete linear convolution operator in $l_p(L)$, with $1\le p\le\infty$ and
$L$ finite. Then,
$\mathcal{K}_{\text{\tiny$\triangle t$}}$
is bounded and
\begin{eqnarray}
\label{th-dis-spe-fin}
\rho(\mathcal{K}_{\text{\tiny$\triangle t$}})
=\rho\left(
\bbb{K}\left( \frac{1}{\triangle t}
\frac{\alpha_{\nu}}{\beta_{\nu}} \right)
\right).
\end{eqnarray}
\end{theorem}

However, the convergence property of the discrete linear
convolution operator $\mathcal{K}_{\text{\tiny$\triangle t$}}$
on infinite time interval is a little bit complicated,
thus, an important lemma is introduced.
\begin{lemma}
\label{le-dis-con-infin}
Assume that the discrete solvability condition (\ref{dis-sol-con}) is satisfied.
Let $b(z)\ne 0$, $\forall$ $|z|=1$.
If $\sp(M_{\mathcal{B}}, -\triangle t M_{\mathcal{A}})\subset\text{\r{S}}$
and $\infty \in S$,
then $\mathcal{K}_{\text{\tiny$\triangle t$}}$ is bounded in $l_p(\infty)$.
Where $S$ represents the absolute stability region of the corresponding linear multistep
formula, and $\text{\r{S}}$ denotes the interior of $S$.
\end{lemma}
\emph{Proof}:
According to the Young's inequality for discrete convolution product \cite{HLP78},
we only need to prove that the kernel $\kappa_{\text{\tiny$\triangle t$}}$
of the discrete convolution operator $\mathcal{K}_{\text{\tiny$\triangle t$}}$
is a $l_1$-sequence.

Denote $\theta^{(-1)}_{\text{\tiny$\triangle t$}}$,
$\theta_{\text{\tiny$\triangle t$}}$ and $\zeta_{\text{\tiny$\triangle t$}}$
as the sequences whose discrete Laplace transforms are
$(a(s)M_{\mathcal{B}}+\triangle t b(s)M_{\mathcal{A}})^{-1}s^k$,
$s^{-k}(a(s)M_{\mathcal{B}}+\triangle t b(s)M_{\mathcal{A}})$,
and $s^{-k}(a(s)N_{\mathcal{B}}+\triangle t b(s)N_{\mathcal{A}})$,
respectively.
Note that the discrete Laplace transform of $\kappa_{\text{\tiny$\triangle t$}}$
satisfies (\ref{dis-lap-mat}), then
\begin{eqnarray*}
\kappa_{\text{\tiny$\triangle t$}} =
\theta^{(-1)}_{\text{\tiny$\triangle t$}} \star
\zeta_{\text{\tiny$\triangle t$}}.
\end{eqnarray*}
Hence, $\kappa_{\text{\tiny$\triangle t$}}$ is a $l_1$-sequence if
$\theta^{(-1)}_{\text{\tiny$\triangle t$}}$ and
$\zeta_{\text{\tiny$\triangle t$}}$ are $l_1$-sequence.
Obviously, $\theta_{\text{\tiny$\triangle t$}}$ and
$\zeta_{\text{\tiny$\triangle t$}}$ are $l_1$-sequence.
Furthermore, according to the Wiener's inversion theorem \cite{Lub83,MieNev87},
$\theta^{(-1)}_{\text{\tiny$\triangle t$}}$ is a
$l_1$-sequence if
\begin{eqnarray}
\label{l1-cond}
\left| a(s)M_{\mathcal{B}}+\triangle t b(s)M_{\mathcal{A}} \right|\ne 0,
\ \forall\ |s|\geq 1.
\end{eqnarray}

Now, we prove the conditions of Lemma \ref{le-dis-con-infin} lead to (\ref{l1-cond}).
Suppose (\ref{l1-cond}) is not satisfied, i.e. there exists a $s_0$ with $|s_0|\ge 1$
such that
\begin{eqnarray}
\label{l1-cond-op}
\left| a(s_0)M_{\mathcal{B}}+\triangle t b(s_0)M_{\mathcal{A}} \right| = 0.
\end{eqnarray}
Considering $M_{\mathcal{B}}$ can be singular or nonsingular, we remark
that: when $b(s_0)\ne 0$, $M_{\mathcal{B}}$ can be either form; when $b(s_0)=0$,
there must be a fact $a(s_0)\ne 0$ due to no common roots for $a(s)$ and $b(s)$,
therefore, $M_{\mathcal{B}}$ must be singular
to make (\ref{l1-cond-op}) satisfied. In order to keep on the discussion, we divide
$b(s_0)$ into two cases, i.e. $b(s_0)\ne 0$ and $b(s_0)= 0$.

If $b(s_0)\ne 0$, (\ref{l1-cond-op}) is equivalent to
\begin{eqnarray}
\label{l1-cond-op-v1}
\left| \frac{a}{b}(s_0)M_{\mathcal{B}}+\triangle t M_{\mathcal{A}} \right| = 0.
\end{eqnarray}
Obviously, (\ref{l1-cond-op-v1}) leads to
\begin{eqnarray}
\label{l1-cond-op-v2}
\frac{a}{b}(s_0)\in\sp(M_{\mathcal{B}}, -\triangle t M_{\mathcal{A}}).
\end{eqnarray}
Since $|s_0|\ge 1$, we have $\frac{a}{b}(s_0)\notin \text{\r{S}}$.
Meanwhile, (\ref{l1-cond-op-v2}) and condition
$\sp(M_{\mathcal{B}}, -\triangle t M_{\mathcal{A}})\subset\text{\r{S}}$
leads to $\frac{a}{b}(s_0)\in \text{\r{S}}$ which contradicts.

If $b(s_0)= 0$, (\ref{l1-cond-op}) is satisfied directly for the singularity of matrix
$M_{\mathcal{B}}$. Moreover, we have $\frac{a}{b}(s_0)=\infty\notin S$ since $|s_0|\ge 1$
which contradicts with the condition $\infty\in S$.
\hfill $\square$

Based on Lemma \ref{dis-spe-inf}, Lemma \ref{le-dis-con-infin}, the
definition of absolute stability region and the maximal principle of
complex function, we can easily obtain the convergence property of
$\mathcal{K}_{\text{\tiny$\triangle t$}}$
on infinite time interval.

\begin{theorem}
\label{th-dis-con-infin}
Assume that the discrete solvability condition (\ref{dis-sol-con}) is satisfied.
Let $b(z)\ne 0$, $\forall$ $|z|=1$.
If $\sp(M_{\mathcal{B}}, -\triangle t M_{\mathcal{A}})\subset\text{\r{S}}$
and $\infty \in S$, consider
$\mathcal{K}_{\text{\tiny$\triangle t$}}$ as a discrete linear convolution operator
in $l_p(\infty)$, with $1\le p\le\infty$. Then
\begin{align}
\rho(\mathcal{K}_{\text{\tiny$\triangle t$}})&=
\sup\{\rho(\bbb{K}(s))\ |\ \triangle ts\in\ceals\backslash\text{\r{S}}\} \label{aaa}\\
&=\sup _{\triangle ts\in\partial S} \rho(\bbb{K}(s)). \label{qqq}
\end{align}
\end{theorem}

\section{Discrete-Time Accelerated Block SOR Method}
\label{sec-dabsor}
The splittings  of  matrices $\mathcal{B}$ and
$\mathcal{A}$ in linear DAEs (\ref{dis-dae-LCC})
are given by
\begin{align}
\mathcal{B} &= M_{\mathcal{B}}-N_{\mathcal{B}}=
\left(\begin{array}{cc}
I&0\\0&0
\end{array}\right) - 0\label{gsor-sin}\\
\intertext{and}
\mathcal{A} &= M_{\mathcal{A}}-N_{\mathcal{A}}=
\left(\begin{array}{cc}
\frac{1}{\omega}A&0\\-B^{T}&\frac{1}{\tau}Q
\end{array}\right) - \left(\begin{array}{cc}
\left(\frac{1}{\omega}-1\right)A&-B\label{gsor-sin-2}\\
0&\frac{1}{\tau}Q
\end{array}\right),
\end{align}
where $Q\in\mathbb{R}^{l\times l}$ is a prescribed
symmetric positive definite matrix as the preconditioner
of the Schur complement matrix $B^TA^{-1}B$. Applying the
generalized successive overrelaxation (GSOR)
technique in \cite{gsor}, we can define the
following discrete-time waveform relaxation method, called the
{\em discrete-time accelerated block successive
overrelaxation} ({\bf DABSOR}) method,
for solving linear DAEs (\ref{dis-dae-LCC}) derived from
time-dependent Stokes equations (\ref{time-stokes}).

\begin{method}
\label{dwr-gsor}
{\sc (The DABSOR Method)}\\
For solving linear constant coefficient DAEs (\ref{dis-dae-LCC})
on time interval $[\txt{T}_1,\txt{T}_2]$, divide the time interval
into $L$ equal distance time steps, and compute the solution of
(\ref{dis-dae-LCC}) on each of the $L$ time levels in
$(\txt{T}_1,\txt{T}_2]$. Let $Q\in\mathbb{R}^{l\times l}$ be
a symmetric positive definite matrix preconditioning the Schur
complement matrix $B^TA^{-1}B$. For two positive integers $r$
and $l$, let $\seq{x}{(0)}$, $\seq{f}{}\in l_p(L;\ceals^r)$ and
$\seq{y}{(0)}$, $\seq{g}{}\in l_p(L;\ceals^l)$ be the initial
iterative vector sequences and the vector sequences derived from
the vector values on each time level of the right hand side of
the linear DAEs (\ref{dis-dae-LCC}). $x_0$, $\ldots$,
$x_{\nu-1}\in\ceals^r$ and $y_0$, $\ldots$, $y_{\nu-1}\in\ceals^l$
are the initial vector values of the iterative vector sequences.
Then:
\begin{itemize}
    \item[For] $k=1,2,\ldots$, untill vector sequences $\seq{x}{(k)}$ and $\seq{y}{(k)}$ converge
        to the exact solution $\seq{x}{}$ and $\seq{y}{}$ of the discrete system derived from
        discretizing the linear DAEs (\ref{dis-dae-LCC}) by linear multistep formulae, compute:
    \item[]
        \begin{itemize}
            \item[For] $n=0:1:L-1$, solve the following linear systems on each time level

                \begin{numcases}{}
                    (\frac{\alpha_{\nu}}{\triangle t}I+\frac{\beta_{\nu}}{\omega}A)x_{n+\nu}^{(k)}=
                    \nonumber\\\qquad\qquad
                    \sum \limits_{j=0}^{\nu}\beta_{j}((\frac{1}{\omega}-1)Ax_{n+j}^{(k-1)}-By_{n+j}^{(k-1)}+f_{n+j})
                    -\sum \limits_{j=0}^{\nu-1}
                    (\frac{\alpha_{j}}{\triangle t}I+\frac{\beta_{j}}{\omega}A)x_{n+j}^{(k)},\nonumber\\
                    \frac{\beta_{\nu}}{\tau}Qy_{n+\nu}^{(k)} = \sum \limits_{j=0}^{\nu}\beta_{j}
                    (B^Tx_{n+j}^{(k)}+\frac{1}{\tau}Qy_{n+j}^{(k-1)}+g_{n+j}) -
                    \sum \limits_{j=0}^{\nu-1} \frac{\beta_{j}}{\tau}Qy_{n+j}^{(k)}.\nonumber
                \end{numcases}
            \item[End]
        \end{itemize}
    \item[End]
\end{itemize}
\end{method}

The DABSOR method can be rewritten into the following matrix form
{\setlength{\multlinegap}{0pt}
\begin{multline}
\label{dabsor-matrix-form}
    \left(\begin{array}{cccccc}
        C_{\nu} & & & & & \\
        C_{\nu-1} & C_{\nu} & & & & \\
        \vdots & \ddots & \ddots & & & \\
        C_0 & \cdots & C_{\nu-1} & C_{\nu} & & \\
         & \ddots & \ddots & \ddots & \ddots & \\
         & & C_0 & \cdots & C_{\nu-1} & C_{\nu} \\
    \end{array}\right)
    \left(\begin{array}{c}
        z_{\nu}^{(k)}\\
        z_{\nu-1}^{(k)}\\
        \vdots\\
        z_{2\nu}^{(k)}\\
        \vdots\\
        z_{L+\nu-1}^{(k)}
    \end{array}\right)=\\
    \left(\begin{array}{cccccc}
        D_{\nu} & & & & & \\
        D_{\nu-1} & D_{\nu} & & & & \\
        \vdots & \ddots & \ddots & & & \\
        D_0 & \cdots & D_{\nu-1} & D_{\nu} & & \\
         & \ddots & \ddots & \ddots & \ddots & \\
         & & D_0 & \cdots & D_{\nu-1} & D_{\nu} \\
    \end{array}\right)
    \left(\begin{array}{c}
        z_{\nu}^{(k-1)}\\
        z_{\nu-1}^{(k-1)}\\
        \vdots\\
        z_{2\nu}^{(k-1)}\\
        \vdots\\
        z_{L+\nu-1}^{(k-1)}
    \end{array}\right)+
    \left(\begin{array}{c}
        \bm{b}_{\nu}\\
        \bm{b}_{\nu-1}\\
        \vdots\\
        \bm{b}_{2\nu}\\
        \vdots\\
        \bm{b}_{L+\nu-1}
    \end{array}\right)
\end{multline}}
with
\begin{eqnarray*}
    C_j = \frac{\alpha_j}{\triangle t} M_{\mathcal{B}} + \beta_j M_{\mathcal{A}}
    \quad\text{ºÍ}\quad D_j = \beta_j N_{\mathcal{A}},\quad j = 0,1,\ldots,\nu,
\end{eqnarray*}
here $M_{\mathcal{B}}$, $M_{\mathcal{A}}$ and $N_{\mathcal{A}}$ are defined in
(\ref{gsor-sin})-(\ref{gsor-sin-2}). Moreover,
$z_{\nu+n}^{(k)}={x_{\nu+n}^{(k)}\choose y_{\nu+n}^{(k)}}$~($n=0,1,\ldots,L-1$),
and
\begin{eqnarray*}
    \bm{b}_{\nu+n} = \left\{
    \begin{array}{l}
        \sum \limits_{j=0}^{\nu} \beta_j b_{n+j}
        + \sum \limits_{j=0}^{\nu-1-n} (D_j-C_j)z_{n+j},\quad \text{if}\quad n=0,1,\ldots,\nu-1,\\
        \\
        \sum \limits_{j=0}^{\nu} \beta_j b_{n+j},\quad \text{if}\quad n=\nu,\nu+1,\ldots,L+\nu-1,
    \end{array}
    \right.
\end{eqnarray*}
with
\begin{eqnarray*}
    b_{\nu+n}={f_{\nu+n}\choose g_{\nu+n}}, n=0,1,\ldots,L-1,\ \text{and}\
    z_{n}={x_{n}\choose y_{n}}, n=0,1,\ldots,\nu-1.
\end{eqnarray*}

\begin{remark}
\label{rk-dabsor-alg}
To be theoretical, the length of the simulation time interval
$[\txt{T}_1,\txt{T}_2]$ in the DABSOR method can be infinite,
i.e. $L$ possibly infinite. In this case, the matrix-vector
multiplications in (\ref{dabsor-matrix-form}) are essentially
the discrete linear convolution between certain matrix and
vector sequences with infinite length. However, in actual
application of the DABSOR method, no computer can deal with
infinite length time interval. Therefore, the length of
$[\txt{T}_1,\txt{T}_2]$ is considered to be finite
in the sequel.
\end{remark}

\subsection{Convergence Domain of Relaxation Parameters}
\label{subsec-dis-con-dom}

The convergence property of the DABSOR method
is described in the following theorem,
which precisely determines the convergence domain
of the DABSOR method with respect to the
relaxation parameters $\omega$ and $\tau$.
For practical application, only the
finite time interval case is studied.
\begin{theorem}
\label{th-dis-con-domain}
Consider the linear DAEs (\ref{dis-dae-LCC}) and the
corresponding DABSOR method, i.e. Method \ref{dwr-gsor},
on finite time interval. Denote the smallest
and the largest eigenvalues of the matrix
$A$ by $\eta_{\min}$ and $\eta_{\max}$, and
those of the matrix $(B^TB)^{-1}Q$ by
$\mu_{\min}$ and $\mu_{\max}$, respectively.
Then the iterative sequence given by the DABSOR
method is convergent, provided
\begin{enumerate}
\item[{\rm (a)}] when
$0\le \frac{\sigma}{\eta_{\max}}
\le\frac{\sigma}{\eta_{\min}}$,
$$
    0<\tau<2\eta_{\min}\mu_{\min}
    \left(\frac{2}{\omega}+
    \frac{\sigma}{\eta_{\max}}-1
    \right),\quad 0<\omega<2;
$$
\item[{\rm (b)}] when
$-1<\frac{\sigma}{\eta_{\min}}
\le\frac{\sigma}{\eta_{\max}}<0$,
$$
    0<\tau<2\eta_{\min}\mu_{\min}
    \left(\frac{2}{\omega}+
    \frac{\sigma}{\eta_{\min}}-1
    \right),
    \quad 0<\omega<
    \frac{2\eta_{\min}}{\eta_{\min}-\sigma};
$$
\item[{\rm (c)}] when
$\frac{\sigma}{\eta_{\min}}
\le\frac{\sigma}{\eta_{\max}}<-1$,
$$
    2\eta_{\min}\mu_{\min}\left(\frac{2}{\omega}+
    \frac{\sigma}{\eta_{\max}}-1
    \right)<\tau<0,\quad
    \frac{2\eta_{\max}}{\eta_{\max}-\sigma}<\omega<2.
$$
\end{enumerate}
Here, $\sigma = \frac{1}{\triangle t}
\frac{\alpha_{\nu}}{\beta_{\nu}}$.
\end{theorem}
\emph{Proof}:
According to Theorem \ref{th-dis-con-fin}, we know that
the spectral radius of the DABSOR method
is given by
\begin{eqnarray*}
\rho(\mathcal{K}_{\text{\tiny$\triangle t$}})
&=&
\rho\left(
\bbb{K}\left( \frac{1}{\triangle t}
\frac{\alpha_{\nu}}{\beta_{\nu}} \right)
\right)=\rho\left( \bbb{K}(\sigma) \right)\\
&=&\rho\left( \left(
\begin{array}{rl}
\sigma I+\frac{1}{\omega}A&0\\
-B^{T}&\frac{1}{\tau}Q
\end{array}
\right)^{-1}
\left(\begin{array}{cc}\left(
\frac{1}{\omega}-1\right)A&-B\\
0&\frac{1}{\tau}Q
\end{array}\right)
\right).
\end{eqnarray*}
Let $\lambda$ be an eigenvalue of the matrix
\begin{eqnarray}
\label{dis-iter-matrix}
\left(\begin{array}{rl}
\sigma I+\frac{1}{\omega}A&0\\-B^{T}&\frac{1}{\tau}Q
\end{array}
\right)^{-1}
\left(\begin{array}{cc}
\left(\frac{1}{\omega}-1\right)A&-B\\
0&\frac{1}{\tau}Q
\end{array}\right)
\end{eqnarray}
and $\left({x\atop y}\right)$ be the
corresponding eigenvector. Then we have
$$
\left(\begin{array}{cc}
\left(\frac{1}{\omega}-1\right)A&-B\\
0&\frac{1}{\tau}Q
\end{array}\right)
\left(\begin{array}{c}x\\y\end{array}\right)=
\lambda \left(\begin{array}{rl}
\sigma I+\frac{1}{\omega}A&0\\-B^{T}&\frac{1}{\tau}Q
\end{array}\right)
\left(\begin{array}{c}x\\y\end{array}\right),
$$
or equivalently
\begin{eqnarray}
\label{dis-ex-eqn-1}
\left\{
    \begin{array}{l}
    \left(1-\omega-\lambda\right)Ax
    -\lambda \omega \sigma x
    =\omega By,\\
    \left(\lambda-1\right)Qy=\lambda\tau B^{T}x.
    \end{array}
    \right.
\end{eqnarray}
Without loss of generality, the vector $x$ is normalized
such that $x^{*}x=1$. Here and in the
sequel, $x^{*}$ is used to denote the conjugate
transpose of the vector $x$. Denote by
$$\gamma_{a}= x^{*}Ax.$$
If $\left(1-\omega-\lambda\right)Ax -
\lambda \omega \sigma x = 0$, then we have
\begin{eqnarray*}
\left(1-\omega-\lambda\right)\gamma_{a} -
\lambda \omega \sigma = 0,
\end{eqnarray*}
or equivalently,
\begin{eqnarray*}
\lambda=\frac{1-\omega}{1+\omega\delta},
\end{eqnarray*}
where $\delta=\frac{\sigma}{\gamma_{a}}$. Thus, we have
\begin{eqnarray*}
\left\{
    \begin{array}{l}
    0=By,\\
    \left(\lambda-1\right)Qy=\lambda\tau B^{T}x.
    \end{array}
    \right.
\end{eqnarray*}
It then follows that $y=0$ and $x\in\zspace(B^T)$, where
$\zspace(B^T)$ represents the null space of the matrix
$B^T$. Hence, $\lambda=\frac{1-\omega}{1+\omega\delta}$
is an eigenvalue of $\bbb{K}(\sigma)$ with
corresponding eigenvector $(x^*,0^*)^*$, where
$x\in\zspace(B^T)$.

Similar to the analysis in \cite{BaiYang12}, we find that
$\lambda=1-\omega$ can also be an
eigenvalue of the matrix in (\ref{dis-iter-matrix}).

Based on the previous cases, two conditions
\begin{eqnarray}
\label{dis-neq1}
|1-\omega|<1\quad\text{and}\quad
\left| \frac{1-\omega}{1+\omega\delta} \right|<1
\end{eqnarray}
should be satisfied to guarantee the convergence of
the DABSOR method.

If $\lambda\neq \frac{1-\omega}{1+\omega\delta}$,
$1-\omega$, then by solving $y$ from the second equation
in (\ref{dis-ex-eqn-1}) we can obtain
$$
y=\frac{\lambda\tau}{\lambda -1}Q^{-1}B^{T}x.
$$
After substituting this equality into the first
equation in (\ref{dis-ex-eqn-1}) we have
$$
\left(1-\omega-\lambda\right)Ax
= \frac{\lambda\omega\tau}{\lambda -1}BQ^{-1}
B^{T}x+\lambda \omega \sigma x.
$$
Premultiplying
$x^{*}$ from left on both sides of the above
equality leads to
\begin{eqnarray}
\label{ex-eqn-2}
\left(1-\omega-\lambda\right)x^{*}Ax
= \frac{\lambda\omega\tau}{\lambda -1}
x^{*}BQ^{-1}B^{T}x+\lambda \omega \sigma.
\end{eqnarray}
Denote by
$$\quad \gamma_{q}=x^{*}BQ^{-1}B^{T}x.$$
Since $x\notin\zspace(B^{T}) $, we obtain
from (\ref{ex-eqn-2}) that
$$
\left(\omega \sigma + \gamma_{a}\right)\lambda^{2}+
\left( \tau\omega \gamma_{q}-\omega \sigma+
\omega \gamma_{a}-2\gamma_{a}
\right)\lambda+\gamma_{a}\left(1-\omega \right)=0,
$$
with the notation
$$
\quad \gamma= \frac{\gamma_{q}}{\gamma_{a}},
$$
the above quadratic polynomial can be rewritten
into the form
\begin{eqnarray}
\label{dis-qua}
\lambda^{2}+\phi\lambda+\psi=0,
\end{eqnarray}
where
\begin{eqnarray*}
\phi &=& \frac{\tau\omega \gamma-
\omega \delta+\omega-2}
{1+\omega \delta}\\
\end{eqnarray*}
and
\begin{eqnarray*}
\psi &=&
\frac{1-\omega}{1+\omega \delta}.
\end{eqnarray*}

Based on the location of zeros of quadratic polynomial
(\ref{dis-qua}) \cite{Mil71} and following the steps
in \cite{BaiYang12}, we obtain
\begin{enumerate}
\item[(a)]
When $ \delta\ge 0$, $\omega$ and $\tau$ satisfy
$$
    0<\tau<\frac{2[2+\omega(\delta-1)]}{\omega\gamma},
    \quad 0< \omega<2.
$$

\item[(b)]
When $ -1<\delta<0$, $\omega$ and $\tau$ satisfy
$$
    0<\tau<\frac{2[2+\omega(\delta-1)]}{\omega\gamma},
    \quad 0<\omega<\frac{2}{1-\delta}.
$$

\item[(c)]
When $\delta<-1$, $\omega$ and $\tau$ satisfy
$$
    \frac{2[2+\omega(\delta-1)]}{\omega\gamma}<\tau<0,
    \quad \frac{2}{1-\delta}<\omega<2.
$$
\end{enumerate}
Recalling that $\delta=\frac{\sigma}{\gamma_{a}}$,
$\gamma=\frac{\gamma_{q}} {\gamma_{a}}$, and
$\gamma_{a}\in [\eta_{\min}, \eta_{\max}]$,
$\gamma_{q} \in [\frac{1}{\mu_{\max}},
\frac{1}{\mu_{\min}}]$, we can easily calculate
the smallest and the largest bounds about
$\delta$ and $\gamma$, denoted by $\delta_{\min}$,
$\gamma_{\min}$ and $\delta_{\max}$,
$\gamma_{\max}$, respectively, as follows:
\begin{enumerate}
\item[(i)] when $\sigma\ge 0$, it holds that
$$
    \begin{array}{l}
    \delta_{\min}=\frac{\sigma}{\eta_{\max}},
    \quad \delta_{\max}= \frac{\sigma}{\eta_{\min}},
    \quad \mbox{and}\quad
    \gamma_{\min}=\frac{1}{\eta_{\max}\mu_{\max}},
    \quad \gamma_{\max}=
    \frac{1}{\eta_{\min}\mu_{\min}};
    \end{array}
$$
\item[(ii)] when $\sigma<0$, it holds that
$$
    \begin{array}{l}
    \delta_{\min}=\frac{\sigma}{\eta_{\min}},
    \quad \delta_{\max}= \frac{\sigma}{\eta_{\max}},
    \quad \mbox{and}\quad
    \gamma_{\min}=\frac{1}{\eta_{\max}\mu_{\max}},
    \quad \gamma_{\max}=
    \frac{1}{\eta_{\min}\mu_{\min}}.
    \end{array}
$$
\end{enumerate}
By making use of these bounds, from the feasible
domain about $\omega$ and $\tau$ determined
in (a)-(c), we can straightforwardly obtain the
following convergence domains for
the DABSOR method:
\begin{enumerate}
\item[(a)]
when $ 0\le \delta_{\min}\le \delta_{\max}$,
$\omega$ and $\tau$ satisfy
$$
    0<\tau<\frac{2[2+\omega(\delta_{\min}-1)]}
    {\omega\gamma_{\max}}, \quad 0< \omega<2;
$$

\item[(b)]
when $ -1<\delta_{\min}\le \delta_{\max}<0$,
$\omega$ and $\tau$ satisfy
$$
    0<\tau<\frac{2[2+\omega(\delta_{\min}-1)]}
    {\omega\gamma_{\max}},
    \quad 0< \omega<\frac{2}{1-\delta_{\min}};
$$

\item[(c)]
when $\delta_{\min}\le\delta_{\max}<-1$,
$\omega$ and $\tau$ satisfy
$$
    \frac{2[2+\omega(\delta_{\max}-1)]}
    {\omega\gamma_{\max}}<\tau<0,
    \quad \frac{2}{1-\delta_{\max}}<\omega<2.
$$
\end{enumerate}
\hfill$\square$

From the proof of Theorem \ref{th-dis-con-domain},
we immediately get the following corollary.
\begin{corollary}
\label{dis-eig-rel}
The nonzero eigenvalues of the matrix $\bbb{K}(\sigma)$
are given by $\lambda = \frac{1-\omega}{1+\omega\delta}$,
$\lambda = 1-\omega$ or
\begin{eqnarray}
\label{dis-qua-eig}
\lambda = \frac{1}{2(1+\omega\delta)}\left[
-(\tau\omega \gamma- \omega \delta+\omega-2)\pm
\sqrt {(\tau\omega \gamma- \omega \delta+\omega-2)^2
-4(1-\omega)(1+\omega\delta) }
\right],
\end{eqnarray}
where $\gamma$ and $\delta$ are the same as in the proof
of Theorem \ref{th-dis-con-domain}.
\end{corollary}

\subsection{The Optimal Iterative Parameters and Convergence Factor}
\label{subsec-dis-opt-par}
In this subsection, the optimal iterative parameters and
the corresponding optimal convergence factor of the
DABSOR method on finite time interval is discussed.
Since the linear multistep formula selected in the
DABSOR method always leads to $\sigma >0$, the condition
$\sigma >0$ is added in the optimality discussion
of the DABSOR method.

Follow the notations in Section \ref{subsec-dis-con-dom},
according Theorem \ref{th-dis-con-domain}, we know
that the iteration parameters $\omega$ and $\tau$ of the
DABSOR method must satisfy
\begin{eqnarray}
\label{dis-interval}
    0<\tau<\frac{2[2+\omega(\delta_{\min}-1)]}
    {\omega\gamma_{\max}}, \quad 0< \omega<2.
\end{eqnarray}

Due to the definition of $\delta$ and the symmetric positive
definite matrix block $A$, the condition $\sigma >0$
leads to $\delta >0$. Therefore, the sequential discussion
is divided into two cases $\delta > 1$ and
$0<\delta\le 1$.

For the case $\delta >1$, we denote the following functions as
\begin{eqnarray*}
\left\{
    \begin{array}{l}
    f_1(\omega,\tau,\gamma,\delta)=
    \frac{1}{2(1+\omega\delta)}\left[
    2-\omega+\omega\delta-\tau\omega \gamma+
    \sqrt {(2-\omega+\omega\delta-\tau\omega \gamma)^2
    -4(1-\omega)(1+\omega\delta) }
    \right],\\
    \qquad\qquad\qquad \text{for}\quad\omega>
    \frac{4\tau\gamma}{(\tau\gamma-\delta+1)^2+4\delta},
    \quad\tau\gamma<\delta-1,\\
    \qquad\qquad\qquad \text{ or}\quad
    \frac{4\tau\gamma}{(\tau\gamma-\delta+1)^2+4\delta}<
    \omega<
    \frac{2}{\tau\gamma-\delta+1},\quad \delta-1<
    \tau\gamma<\delta+1;\\
    f_2(\omega,\tau,\gamma,\delta)=
    \frac{1}{2(1+\omega\delta)}\left[
    \tau\omega \gamma- \omega \delta+\omega-2 +
    \sqrt {(\tau\omega \gamma- \omega \delta+\omega-2)^2
    -4(1-\omega)(1+\omega\delta) }
    \right],\\
    \qquad\qquad\qquad \text{for}\quad\omega>
    \frac{4\tau\gamma}{(\tau\gamma-\delta+1)^2+4\delta},
    \quad\tau\gamma>\delta+1,\\
    \qquad\qquad\qquad \text{ or}\quad\omega>
    \frac{2}{\tau\gamma-\delta+1},\quad \delta-1<
    \tau\gamma<\delta+1;\\
    f_3(\omega,\tau,\gamma,\delta)=g(\omega,\delta)=
    \sqrt {\frac{1-\omega}{1+\omega\delta}},\\
    \qquad\qquad\qquad \text{for}\quad\omega<
    \frac{4\tau\gamma}{(\tau\gamma-\delta+1)^2+4\delta}.
    \end{array}
    \right.
\end{eqnarray*}
The above functions are induced from calculating the
magnitudes of the nonzero eigenvalues of the
matrix $\bbb{K}(\sigma)$ given in the Corollary \ref{dis-eig-rel}
based on the following three cases:
\begin{enumerate}
\item[(i)] $2-\omega+\omega\delta-\tau\omega \gamma>0$,
$(2-\omega+\omega\delta-\tau\omega \gamma)^2
-4(1-\omega)(1+\omega\delta)>0$;
\item[(ii)] $2-\omega+\omega\delta-\tau\omega \gamma<0$,
$(2-\omega+\omega\delta-\tau\omega \gamma)^2
-4(1-\omega)(1+\omega\delta)>0$;
\item[(iii)] $(2-\omega+\omega\delta-\tau\omega \gamma)^2
-4(1-\omega)(1+\omega\delta)\le 0$.
\end{enumerate}
The first two cases correspond to the positive discriminant
of the quadratic polynomial (\ref{dis-qua}) and the sign
of the term $2-\omega+\omega\delta-\tau\omega \gamma$.
Meanwhile, the third case corresponds to the
non-positive discriminant. Further investigating each of
these cases, together with the intervals given in
(\ref{dis-interval}) with respect to $\omega$ and $\tau$,
leads to the definitions of functions $f_j(\omega,\tau,\gamma,\delta)$
($j=1,2,3$) respectively.
\begin{assumption}
\label{dis-assumption-1}
$\frac{\omega\sigma}{\tau}\notin\sp(Q^{-1}B^TB)$.
\end{assumption}

On account of Assumption \ref{dis-assumption-1}, similar to the analysis in
\cite{BaiYang12}, we have
\begin{eqnarray*}
(1-\omega)\notin\sp(\bbb{K}(\sigma)).
\end{eqnarray*}
At the same time, the restrictions on $\omega$ and $\tau$ in
the definitions of functions $f_j(\omega,\tau,\gamma,\delta)$
($j=1,2,3$) make it holds that
\begin{eqnarray*}
f_j(\omega,\tau,\gamma,\delta)\ge
\sqrt {\frac{1-\omega}{1+\omega\delta}}
\ge \frac{1-\omega}{1+\omega\delta},
\quad j=1,2,3.
\end{eqnarray*}

Now, we discuss the monotone properties of the functions
$f_j(\omega,\tau,\gamma,\delta)$ ($j=1,2,3$) with respect
to $\gamma$ and $\omega$.

By specific computation, we get
\begin{eqnarray*}
\left\{
    \begin{array}{l}
    \frac{\partial f_1(\omega,\tau,\gamma,\delta)}
    {\partial \gamma}=
    -\frac{\tau\omega}{2(1+\omega\delta)}\left[
    1+\frac{1}{\sqrt {(2-\omega+\omega\delta-\tau\omega \gamma)^2
    -4(1-\omega)(1+\omega\delta)} }
    \right],\\
    \qquad\qquad\qquad \text{for}\quad\omega>
    \frac{4\tau\gamma}{(\tau\gamma-\delta+1)^2+4\delta},
    \quad\tau\gamma<\delta-1,\\
    \qquad\qquad\qquad \text{ or}\quad
    \frac{4\tau\gamma}{(\tau\gamma-\delta+1)^2+4\delta}<
    \omega<
    \frac{2}{\tau\gamma-\delta+1},\quad \delta-1<
    \tau\gamma<\delta+1;\\
    \frac{\partial f_2(\omega,\tau,\gamma,\delta)}
    {\partial \gamma}=
    \frac{\tau\omega}{2(1+\omega\delta)}\left[
    1+\frac{1}{
    \sqrt {(\tau\omega \gamma- \omega \delta+\omega-2)^2
    -4(1-\omega)(1+\omega\delta)} }
    \right],\\
    \qquad\qquad\qquad \text{for}\quad\omega>
    \frac{4\tau\gamma}{(\tau\gamma-\delta+1)^2+4\delta},
    \quad\tau\gamma>\delta+1,\\
    \qquad\qquad\qquad \text{ or}\quad\omega>
    \frac{2}{\tau\gamma-\delta+1},\quad \delta-1<
    \tau\gamma<\delta+1;\\
    \frac{\partial f_3(\omega,\tau,\gamma,\delta)}
    {\partial \gamma}=0,\\
    \qquad\qquad\qquad \text{for}\quad\omega<
    \frac{4\tau\gamma}{(\tau\gamma-\delta+1)^2+4\delta}.
    \end{array}
    \right.
\end{eqnarray*}
Thus
\begin{eqnarray}
\label{dis-mon-gamma}
\left\{
    \begin{array}{l}
    \frac{\partial f_1(\omega,\tau,\gamma,\delta)}
    {\partial \gamma}<0,\quad \text{for}\quad\omega>
    \frac{4\tau\gamma}{(\tau\gamma-\delta+1)^2+4\delta},
    \quad\tau\gamma<\delta-1,\\
    \qquad\qquad\qquad\quad\ \ \text{ or}\quad
    \frac{4\tau\gamma}{(\tau\gamma-\delta+1)^2+4\delta}<
    \omega<
    \frac{2}{\tau\gamma-\delta+1},\quad \delta-1<
    \tau\gamma<\delta+1;\\
    \frac{\partial f_2(\omega,\tau,\gamma,\delta)}
    {\partial \gamma}>0,
    \quad \text{for}\quad\omega>
    \frac{4\tau\gamma}{(\tau\gamma-\delta+1)^2+4\delta},
    \quad\tau\gamma>\delta+1,\\
    \qquad\qquad\qquad\quad\ \ \text{ or}\quad\omega>
    \frac{2}{\tau\gamma-\delta+1},\quad \delta-1<
    \tau\gamma<\delta+1;\\
    \frac{\partial f_3(\omega,\tau,\gamma,\delta)}
    {\partial \gamma}=0, \quad \text{for}\quad\omega<
    \frac{4\tau\gamma}{(\tau\gamma-\delta+1)^2+4\delta}.
    \end{array}
    \right.
\end{eqnarray}
Based on (\ref{dis-mon-gamma}), we find that:
$f_1(\omega,\tau,\gamma,\delta)$ decreases with
respect to $\gamma$ when
$\omega> \frac{4\tau\gamma}{(\tau\gamma-\delta+1)^2+4\delta}$
and $\tau\gamma<\delta-1$, or when
$\frac{4\tau\gamma}{(\tau\gamma-\delta+1)^2+4\delta}< \omega<
\frac{2}{\tau\gamma-\delta+1}$ and
$\delta-1<\tau\gamma<\delta+1$;
$f_2(\omega,\tau,\gamma,\delta)$ increases with
respect to $\gamma$ when
$\omega>\frac{4\tau\gamma}{(\tau\gamma-\delta+1)^2+4\delta}$ and
$\tau\gamma>\delta+1$, or when
$\omega>\frac{2}{\tau\gamma-\delta+1}$ and
$\delta-1<\tau\gamma<\delta+1$;
$f_3(\omega,\tau,\gamma,\delta)$ is not related to $\gamma$.

Denote $f_j^{nume}$ $(j=1,2,3)$ as the numerator of the
functions $f_j$ $(j=1,2,3)$. Then
\begin{eqnarray*}
\left\{
    \begin{array}{l}
    \frac{\partial f_1^{nume}(\omega,\tau,\gamma,\delta)}
    {\partial \omega}=
    -\left[\tau\gamma-\delta+1
    +\frac{(2-\omega+\omega\delta-\tau\omega \gamma)
    (\tau\gamma-\delta+1)+2(\delta-1-2\omega\delta)}
    {\sqrt {(2-\omega+\omega\delta-\tau\omega \gamma)^2
    -4(1-\omega)(1+\omega\delta)} }
    \right],\\
    \qquad\qquad\qquad \text{for}\quad\omega>
    \frac{4\tau\gamma}{(\tau\gamma-\delta+1)^2+4\delta},
    \quad\tau\gamma<\delta-1,\\
    \qquad\qquad\qquad \text{ or}\quad
    \frac{4\tau\gamma}{(\tau\gamma-\delta+1)^2+4\delta}<
    \omega<
    \frac{2}{\tau\gamma-\delta+1},\quad \delta-1<
    \tau\gamma<\delta+1;\\
    \frac{\partial f_2^{nume}(\omega,\tau,\gamma,\delta)}
    {\partial \omega}=\tau\gamma-\delta+1
    +\frac{(\tau\omega \gamma- \omega \delta+\omega-2)
    (\tau\gamma-\delta+1)-2(\delta-1-2\omega\delta)}
    {\sqrt {(\tau\omega \gamma- \omega \delta+\omega-2)^2
    -4(1-\omega)(1+\omega\delta)} } ,\\
    \qquad\qquad\qquad \text{for}\quad\omega>
    \frac{4\tau\gamma}{(\tau\gamma-\delta+1)^2+4\delta},
    \quad\tau\gamma>\delta+1,\\
    \qquad\qquad\qquad \text{ or}\quad\omega>
    \frac{2}{\tau\gamma-\delta+1},\quad \delta-1<
    \tau\gamma<\delta+1;\\
    \end{array}
    \right.
\end{eqnarray*}
Therefore
\begin{eqnarray*}
\left\{
    \begin{array}{l}
    \frac{\partial f_1(\omega,\tau,\gamma,\delta)}
    {\partial \omega}=
    \frac{(1+\omega\delta)\frac{\partial}{\partial \omega}
    f_1^{nume}-\delta f_1^{nume}}
    {2(1+\omega\delta)^2},\\
    \qquad\qquad\qquad \text{for}\quad\omega>
    \frac{4\tau\gamma}{(\tau\gamma-\delta+1)^2+4\delta},
    \quad\tau\gamma<\delta-1,\\
    \qquad\qquad\qquad \text{ or}\quad
    \frac{4\tau\gamma}{(\tau\gamma-\delta+1)^2+4\delta}<
    \omega<
    \frac{2}{\tau\gamma-\delta+1},\quad \delta-1<
    \tau\gamma<\delta+1;\\
    \frac{\partial f_2(\omega,\tau,\gamma,\delta)}
    {\partial \omega}=
    \frac{(1+\omega\delta)\frac{\partial}{\partial \omega}
    f_2^{nume}-\delta f_2^{nume}}
    {2(1+\omega\delta)^2},\\
    \qquad\qquad\qquad \text{for}\quad\omega>
    \frac{4\tau\gamma}{(\tau\gamma-\delta+1)^2+4\delta},
    \quad\tau\gamma>\delta+1,\\
    \qquad\qquad\qquad \text{ or}\quad\omega>
    \frac{2}{\tau\gamma-\delta+1},\quad \delta-1<
    \tau\gamma<\delta+1;\\
    \frac{\partial f_3(\omega,\tau,\gamma,\delta)}
    {\partial \omega}=\frac{-1}
    {\sqrt {1-\omega}\left( \sqrt {1+\omega\delta} \right)^3},\\
    \qquad\qquad\qquad \text{for}\quad\omega<
    \frac{4\tau\gamma}{(\tau\gamma-\delta+1)^2+4\delta}.
    \end{array}
    \right.
\end{eqnarray*}
Denote
\begin{eqnarray*}
\vartheta(\omega)=\delta\left[ (\eta^2+4\delta)\eta -
\delta(\delta-1)^2 \right]\omega^2 - 2\delta(\eta+\tau\gamma)
\tau\gamma\omega-4\tau^2\gamma^2,
\end{eqnarray*}
where $\eta=\tau\gamma-\delta+1$.
\begin{assumption}
\label{dis-assumption-2}
$\vartheta(\omega)\le 0$ when
$\omega> \frac{4\tau\gamma}{(\tau\gamma-\delta+1)^2+4\delta}$
and $\tau\gamma<\delta-1$, or when
$\frac{4\tau\gamma}{(\tau\gamma-\delta+1)^2+4\delta}< \omega<
\frac{2}{\tau\gamma-\delta+1}$ and
$\delta-1<\tau\gamma<\delta+1$.
\end{assumption}
Then
\begin{eqnarray}
\label{dis-mon-omega}
\left\{
    \begin{array}{l}
    \frac{\partial f_1(\omega,\tau,\gamma,\delta)}
    {\partial \omega}>0,\quad
    \text{when Assumption \ref{dis-assumption-2} is true, and}\\
    \qquad\qquad\qquad\quad\ \ \text{for}\quad\omega>
    \frac{4\tau\gamma}{(\tau\gamma-\delta+1)^2+4\delta},
    \quad\tau\gamma<\delta-1,\\
    \qquad\qquad\qquad\quad\ \ \text{ or}\quad
    \frac{4\tau\gamma}{(\tau\gamma-\delta+1)^2+4\delta}<
    \omega<
    \frac{2}{\tau\gamma-\delta+1},\quad \delta-1<
    \tau\gamma<\delta+1;\\
    \frac{\partial f_2(\omega,\tau,\gamma,\delta)}
    {\partial \omega}>0,
    \quad \text{for}\quad\omega>
    \frac{4\tau\gamma}{(\tau\gamma-\delta+1)^2+4\delta},
    \quad\tau\gamma>\delta+1,\\
    \qquad\qquad\qquad\quad\ \ \text{ or}\quad\omega>
    \frac{2}{\tau\gamma-\delta+1},\quad \delta-1<
    \tau\gamma<\delta+1;\\
    \frac{\partial f_3(\omega,\tau,\gamma,\delta)}
    {\partial \omega}<0, \quad \text{for}\quad\omega<
    \frac{4\tau\gamma}{(\tau\gamma-\delta+1)^2+4\delta}.
    \end{array}
    \right.
\end{eqnarray}
In view of (\ref{dis-mon-omega}), we conclude that:
$f_1(\omega,\tau,\gamma,\delta)$ increases with
respect to $\omega$ when
$\omega> \frac{4\tau\gamma}{(\tau\gamma-\delta+1)^2+4\delta}$
and $\tau\gamma<\delta-1$, or when
$\frac{4\tau\gamma}{(\tau\gamma-\delta+1)^2+4\delta}< \omega<
\frac{2}{\tau\gamma-\delta+1}$ and
$\delta-1<\tau\gamma<\delta+1$;
$f_2(\omega,\tau,\gamma,\delta)$ increases with
respect to $\omega$ when
$\omega>\frac{4\tau\gamma}{(\tau\gamma-\delta+1)^2+4\delta}$ and
$\tau\gamma>\delta+1$, or when
$\omega>\frac{2}{\tau\gamma-\delta+1}$ and
$\delta-1<\tau\gamma<\delta+1$;
$f_3(\omega,\tau,\gamma,\delta)$ decreases with
respect to $\omega$ when
$\omega<\frac{4\tau\gamma}{(\tau\gamma-\delta+1)^2+4\delta}$.

Moreover, when $f_j(\omega,\tau,\gamma,\delta)$ ($j=1,2,3$) are well
defined for positive reals $\gamma$ and $\tau$, we have
\begin{eqnarray}
\label{dis-intersection-1}
\left\{
    \begin{array}{l}
    f_1(\omega,\tau,\gamma,\delta)=f_3(\omega,\tau,\gamma,\delta)
    \quad\text{if} \quad \omega=\frac{4\tau\gamma}
    {(\tau\gamma-\delta+1)^2+4\delta},\\
    f_2(\omega,\tau,\gamma,\delta)=f_3(\omega,\tau,\gamma,\delta)
    \quad\text{if} \quad \omega=\frac{4\tau\gamma}
    {(\tau\gamma-\delta+1)^2+4\delta},\\
    f_1(\omega,\tau,\gamma,\delta)=f_2(\omega,\tau,\gamma,\delta)
    \quad\text{if} \quad \omega=\frac{2}
    {\tau\gamma-\delta+1}.\\
    \end{array}
    \right.
\end{eqnarray}
Denote $\omega(\tau,\gamma)=\frac{4\tau\gamma}
{(\tau\gamma-\delta+1)^2+4\delta}$, then $\omega(\tau,\gamma)$
increases with respect to $\gamma$ when $\tau\gamma<\delta+1$,
and $\omega(\tau,\gamma)$ decreases with respect to
$\gamma$ when $\tau\gamma>\delta+1$.

For two different reals
$\gamma_1$ and $\gamma_2$, we have
\begin{eqnarray}
\label{dis-intersection-2}
f_1(\omega,\tau,\gamma_1,\delta)=f_2(\omega,\tau,\gamma_2,\delta)
\quad\text{if} \quad \omega=\frac{4}
{\tau(\gamma_1+\gamma_2)-2(\delta-1)}.
\end{eqnarray}

In addition, we define the functions
$$
\omega_-(\tau) = \omega(\gamma_{\min},\delta),\quad
\omega_+(\tau) = \omega(\gamma_{\max},\delta),\quad
\omega_0(\tau) = \frac{4}
{\tau(\gamma_{\min}+\gamma_{\max})-2(\delta-1)}.
$$

By applying the Corollary \ref{dis-eig-rel}, after concrete
computations we know that the magnitudes of the nonzero
eigenvalues $\lambda$ of the matrix $\bbb{K}(\sigma)$
can be expressed as following:

when $\tau\gamma<\delta-1$, $|\lambda|=\left| \frac{1-\omega}
{1+\omega\delta} \right|$ or
\begin{eqnarray}
\label{dis-case-1}
|\lambda|=\left\{
    \begin{array}{l}
    f_1(\omega,\tau,\gamma,\delta),\quad\text{for}\quad
    \omega>\frac{4\tau\gamma}
    {(\tau\gamma-\delta+1)^2+4\delta},\\
    g(\omega,\delta),\quad\text{for}\quad
    \omega<\frac{4\tau\gamma}
    {(\tau\gamma-\delta+1)^2+4\delta};
    \end{array}
    \right.
\end{eqnarray}

when $\delta-1<\tau\gamma<\delta+1$, $|\lambda|=\left| \frac{1-\omega}
{1+\omega\delta} \right|$ or
\begin{eqnarray}
\label{dis-case-2}
|\lambda|=\left\{
    \begin{array}{l}
    f_1(\omega,\tau,\gamma,\delta),\quad\text{for}\quad
    \frac{4\tau\gamma}
    {(\tau\gamma-\delta+1)^2+4\delta}<
    \omega<\frac{2}{\tau\gamma-\delta+1},\\
    f_2(\omega,\tau,\gamma,\delta),\quad\text{for}\quad
    \omega>\frac{2}{\tau\gamma-\delta+1},\\
    g(\omega,\delta),\quad\text{for}\quad
    \omega<\frac{4\tau\gamma}
    {(\tau\gamma-\delta+1)^2+4\delta};
    \end{array}
    \right.
\end{eqnarray}

when $\tau\gamma>\delta+1$, $|\lambda|=\left| \frac{1-\omega}
{1+\omega\delta} \right|$ or
\begin{eqnarray}
\label{dis-case-3}
|\lambda|=\left\{
    \begin{array}{l}
    f_2(\omega,\tau,\gamma,\delta),\quad\text{for}\quad
    \omega>\frac{4\tau\gamma}
    {(\tau\gamma-\delta+1)^2+4\delta},\\
    g(\omega,\delta),\quad\text{for}\quad
    \omega<\frac{4\tau\gamma}
    {(\tau\gamma-\delta+1)^2+4\delta}.
    \end{array}
    \right.
\end{eqnarray}

By observing (\ref{dis-case-1})-(\ref{dis-case-2}),
we find that in order
to compute the spectral radius of $\rho(\bbb{K}(\sigma))$
we have to discuss in
the following three cases with respect to the parameter $\tau$:
\begin{enumerate}
\item[(a)] $\tau\le\frac{\delta-1}{\gamma_{\max}}$;
\item[(b)] $\tau\ge\frac{\delta+1}{\gamma_{\min}}$;
\item[(c)] $\frac{\delta-1}{\gamma_{\max}}<\tau
<\frac{\delta+1}{\gamma_{\min}}$.
\end{enumerate}

For the case $0<\delta\le 1$, a similar discussion can be stated
as long as Assumption \ref{dis-assumption-2} is replaced
as the following one.
\begin{assumption}
\label{dis-assumption-3}
$\vartheta(\omega)\le 0$ when
$\frac{4\tau\gamma}{(\tau\gamma-\delta+1)^2+4\delta}< \omega<
\frac{2}{\tau\gamma-\delta+1}$ and $0<\tau\gamma<\delta+1$.
\end{assumption}

Based on the above analysis, we can give a specific demonstration
of the optimal iterative parameters and the corresponding optimal
convergence factor of the DABSOR method.

\begin{theorem}
\label{th-dis-opt}
Consider the linear DAEs (\ref{dis-dae-LCC}) and the
corresponding DABSOR method, i.e. Method \ref{dwr-gsor}, on finite
time interval. Let $\sigma>0$, and take the same notations
in Theorem \ref{th-dis-con-domain}.

For the case $\delta>1$,
let $\hat \tau=\frac{\delta+1}{\sqrt {\gamma_{\min}\gamma_{\max}}}$,
then
$\omega_-(\hat \tau)=\omega_+(\hat \tau)
=\omega_0(\hat \tau),$
and
\begin{eqnarray*}
\left\{
    \begin{array}{l}
    \omega_-(\tau)\le\omega_+(\tau)\le\omega_0(\tau),
    \quad\frac{\delta-1}{\gamma_{\max}}<\tau<\hat \tau,\\
    \omega_0(\tau)\le\omega_+(\tau)\le\omega_-(\tau),
    \quad\hat \tau<\tau<\frac{\delta+1}{\gamma_{\min}}.
    \end{array}
    \right.
\end{eqnarray*}
When $\frac{\delta-1}{\gamma_{\max}}<\tau<\hat \tau$, let
$\gamma_1$, $\gamma_2\in (\gamma_{\min},\gamma_{\max})$
be positive reals satisfying
\begin{eqnarray*}
\left\{
    \begin{array}{l}
    \gamma_1=\sup\left\{
        \gamma\ |\ \tau\gamma<\delta-1\quad\text{and}\quad
        \frac{\delta-1}{\gamma_{\max}}<\tau<\hat \tau\right\},\\
    \gamma_2=\sup\left\{
        \gamma\ |\ \tau\gamma<\delta+1\quad\text{and}\quad
        \frac{\delta-1}{\gamma_{\max}}<\tau<\hat \tau\right\}.
    \end{array}
    \right.
\end{eqnarray*}
When $\hat \tau<\tau<\frac{\delta+1}{\gamma_{\min}}$, let
$\gamma_1$, $\gamma_2\in (\gamma_{\min},\gamma_{\max})$
be positive reals satisfying
\begin{eqnarray*}
\left\{
    \begin{array}{l}
    \gamma_1=\sup\left\{
        \gamma\ |\ \tau\gamma<\delta-1\quad\text{and}\quad
        \hat \tau<\tau<\frac{\delta+1}{\gamma_{\min}}\right\},\\
    \gamma_2=\sup\left\{
        \gamma\ |\ \tau\gamma<\delta+1\quad\text{and}\quad
        \hat \tau<\tau<\frac{\delta+1}{\gamma_{\min}}\right\}.
    \end{array}
    \right.
\end{eqnarray*}
Denote $\hat \omega_0(\tau)=\frac{4}{\tau(\gamma_1+\gamma_2)-2(\delta-1)}$
as the point of intersection of $f_1(\omega,\tau,\gamma_1,\delta)$ and
$f_2(\omega,\tau,\gamma_2,\delta)$,
$\hat \omega_-(\tau) = \omega(\tau,\gamma_1)$ as the point of intersection of
$f_1(\omega,\tau,\gamma_1,\delta)$ and $g(\omega,\delta)$,
$\hat \omega_+(\tau) = \omega(\tau,\gamma_2)$ as the point of intersection of
$f_2(\omega,\tau,\gamma_2,\delta)$ and $g(\omega,\delta)$.
If Assumptions
\ref{dis-assumption-1} and \ref{dis-assumption-2} are
satisfied, then
\begin{enumerate}
\item[(a)] when $\frac{\delta-1}{\gamma_{\max}}<\tau<\frac{\delta+1}
{\sqrt {\gamma_{\min}\gamma_{\max}}}$,
\begin{enumerate}
\item[(i)] for
$\hat \omega_-(\tau)<\hat \omega_+(\tau)<\hat \omega_0(\tau)$,
\begin{eqnarray*}
\rho(\bbb{K}(\sigma))=\left\{
    \begin{array}{l}
    g(\omega,\delta),\quad\text{for}\quad
    0<\omega<\omega_-(\tau),\\
    f_1(\omega,\tau,\gamma_{\min},\delta),\quad\text{for}\quad
    \omega_-(\tau)<\omega<\omega_0(\tau),\\
    f_2(\omega,\tau,\gamma_{\max},\delta),\quad\text{for}\quad
    \omega_0(\tau)<\omega<2;\\
    \end{array}
    \right.
\end{eqnarray*}
\item[(ii)] for
$\hat \omega_0(\tau)<\hat \omega_-(\tau)<\hat \omega_+(\tau)$,
\begin{eqnarray*}
\rho(\bbb{K}(\sigma))=\left\{
    \begin{array}{l}
    g(\omega,\delta),\quad\text{for}\quad
    0<\omega<\omega_-(\tau),\\
    f_1(\omega,\tau,\gamma_{\min},\delta),\quad\text{for}\quad
    \omega_-(\tau)<\omega<2;\\
    \end{array}
    \right.
\end{eqnarray*}
\end{enumerate}
\item[(b)] when $\frac{\delta+1}{\sqrt {\gamma_{\min}\gamma_{\max}}}
<\tau<\frac{\delta+1}{\gamma_{\min}}$,
\begin{eqnarray*}
\rho(\bbb{K}(\sigma))=\left\{
    \begin{array}{l}
    g(\omega,\delta),\quad\text{for}\quad
    0<\omega<\omega_+(\tau),\\
    f_2(\omega,\tau,\gamma_{\max},\delta),\quad\text{for}\quad
    \omega_+(\tau)<\omega<2.\\
    \end{array}
    \right.
\end{eqnarray*}
\end{enumerate}

For the case $0<\delta\le 1$, if Assumptions \ref{dis-assumption-1} and
\ref{dis-assumption-3} are satisfied, then
\begin{enumerate}
\item[(a)] when $\frac{\delta+1}{\gamma_{\max}}<\tau<\frac{\delta+1}
{\sqrt {\gamma_{\min}\gamma_{\max}}}$,
\begin{eqnarray*}
\rho(\bbb{K}(\sigma))=\left\{
    \begin{array}{l}
    g(\omega,\delta),\quad\text{for}\quad
    0<\omega<\omega_-(\tau),\\
    f_1(\omega,\tau,\gamma_{\min},\delta),\quad\text{for}\quad
    \omega_-(\tau)<\omega<\omega_0(\tau),\\
    f_2(\omega,\tau,\gamma_{\max},\delta),\quad\text{for}\quad
    \omega_0(\tau)<\omega<2;\\
    \end{array}
    \right.
\end{eqnarray*}
\item[(b)] when $\frac{\delta+1}{\sqrt {\gamma_{\min}\gamma_{\max}}}
<\tau<\frac{\delta+1}{\gamma_{\min}}$,
\begin{eqnarray*}
\rho(\bbb{K}(\sigma))=\left\{
    \begin{array}{l}
    g(\omega,\delta),\quad\text{for}\quad
    0<\omega<\omega_+(\tau),\\
    f_2(\omega,\tau,\gamma_{\max},\delta),\quad\text{for}\quad
    \omega_+(\tau)<\omega<2.\\
    \end{array}
    \right.
\end{eqnarray*}
\end{enumerate}

Furthermore, for any $\delta>0$, the optimal iterative
parameters $\tau_{\opt}$ and $\omega_{\opt}$ are given by
$$
\tau_{\opt}=\frac{\delta+1}{\sqrt {\gamma_{\min}\gamma_{\max}}}\quad
\text{and}\quad
\omega_{\opt}=\frac{4\sqrt {\gamma_{\min}\gamma_{\max}}}
{(\delta+1)(\gamma_{\min}+\gamma_{\max})-2(\delta-1)
\sqrt {\gamma_{\min}\gamma_{\max}}},
$$
and the corresponding optimal convergence factor of the DABSOR method
is given by
$$
\rho(\bbb{K}(\sigma))_{\opt}=
\frac{\sqrt {\gamma_{\max}}-\sqrt {\gamma_{\min}}}
{\sqrt {\gamma_{\max}}+\sqrt {\gamma_{\min}}},
$$
here, $\delta\in(\delta_{\min},\delta_{\max})$.
\end{theorem}

\begin{remark}
\label{rk-opt-con-curve}
According to the expressions of the optimal iterative parameters
$\tau_{\opt}$ and $\omega_{\opt}$ in Theorem \ref{th-dis-opt},
$(\tau_{\opt},\omega_{\opt})$ is not a single point
in $\omega$-$\tau$ plane, which means that the optimal convergence
factor $\rho(\bbb{K}(\sigma))_{\opt}$ is obtained on
a parameterized curve $(\tau_{\opt}(\delta),\omega_{\opt}(\delta))$
with respect to $\delta\in(\delta_{\min},\delta_{\max})$.
Thus, the above parameterized curve which shows all the optimal
parameters is called {\bf optimal convergence curve}. Moreover,
the properties of the optimal convergence curve are closely
related to the linear multistep formulae selected and the properties
of the matrix block $A$, $B$ and the preconditioner $Q$.
\end{remark}

\subsection{The DABSOR Method with Windowing Technique}
\label{subsec-dis-dwr-wintec}

It is known that there is a typical phenomenon of
the waveform relaxation methods based on matrix splitting
that standard matrix splitting iterative methods for solving linear
algebraic system may not have, that is, during the iterative procedure of the
waveform relaxation methods, the intermediate solutions contain
spurious oscillations with growth of the error and translation
of the oscillating region.

To be specific, according to Theorems \ref{th-dis-con-fin}
and \ref{th-dis-con-infin}, the spectral radius of
$\mathcal{K}_{\text{\tiny$\triangle t$}}$ as a discrete linear convolution
operator on finite time interval is smaller than that on infinite
time interval. Thus, it is reasonable that the waveform relaxation
methods are convergent on finite time interval, and
divergent on infinite time interval. In this case, the iterative
procedure on a sufficient long time interval firstly seems
to diverge, i.e. oscillations appear in large part of
the whole computation time interval. Eventually, the iterative
procedure surely starts to converge, i.e. the length of
time interval with small error extends slowly as the iteration
proceeds. Therefore, the asymptotic convergence behavior is
dictated by Theorem \ref{th-dis-con-fin}. Nevertheless, it takes
a large number of iterative steps to make the region of divergent
behavior receding backward, which predicts a rapid increase to
the computation load.

In order to get around the above shortcoming during long time
interval simulation, an acceleration technique, called {\bf windowing},
is introduced to the waveform relaxation methods.
In fact, windowing is a technique to divide the whole long time
interval into a number of short time subintervals based on certain rules, and
apply the corresponding waveform relaxation methods on each subinterval.
Since the subintervals are short, the number of iterative steps of
the waveform relaxation methods on each
subinterval is smaller than that on the
whole long time interval.
Furthermore, the sum of computation loads on all of the subintervals
is certainly smaller than the computation load while simulating on
the whole long time interval.
To improve computing efficiency,
the DABSOR method is integrated with windowing technique.
\begin{method}
\label{dwr-gsor-win}
{\sc (The DABSOR Method with Windowing Technique)}\\
For solving linear constant coefficient DAEs (\ref{dis-dae-LCC})
on time interval $[\txt{T}_1,\txt{T}_2]$, divide the time interval
into $L$ equal distance time steps, and compute the solution of
(\ref{dis-dae-LCC}) on each of the $L$ time levels in
$\Omega_t=(\txt{T}_1,\txt{T}_2]$. Choosing $N+1$ time levels
$\txt{T}_1=\txt{t}_0<\txt{t}_1<\cdots<\txt{t}_N=\txt{T}_2$ to
divide time interval $\Omega_t$ into $N$ smaller subintervals
$\Omega_t^{(i)}=(\txt{t}_{i-1},\txt{t}_{i}]$, $i=1,2,\ldots,N$,
with $L_i$ time levels in each subinterval $\Omega_t^{(i)}$, and
$\sum _{i=1}^{N} L_i = L$. Let $Q\in\mathbb{R}^{l\times l}$ be
a symmetric positive definite matrix preconditioning the Schur
complement matrix $B^TA^{-1}B$. For two positive integers $r$ and
$l$, let $\seqi{x}{i}{(0)}$, $\seqi{f}{i}{}\in l_p(L_i;\ceals^r)$
and $\seqi{y}{i}{(0)}$, $\seqi{g}{i}{}\in l_p(L_i;\ceals^l)$
be the initial iterative vector sequences and the vector sequences
derived from the vector values on each of the corresponding $L_i$
time levels of the right hand side of the linear DAEs (\ref{dis-dae-LCC}).
$x_0$, $\ldots$, $x_{\nu-1}\in\ceals^r$ and $y_0$, $\ldots$,
$y_{\nu-1}\in\ceals^l$ are the initial vector values of the
iterative vector sequences.
Then:
\begin{itemize}
    \item[For] $i=1,2,\ldots,N$, on each subinterval $\Omega_t^{(i)}$, compute
        \begin{itemize}
            \item[For] $k=1,2,\ldots$£¬untill vector sequences $\seqi{x}{i}{(k)}$ and $\seqi{y}{i}{(k)}$
            converge to the exact solution $\seqi{x}{i}{}$ and $\seqi{y}{i}{}$ of the discrete system
            derived from discretizing the linear DAEs (\ref{dis-dae-LCC}) by linear multistep formulae,
            compute
            \item[]
                \begin{itemize}
                    \item[For] $n=\sum _{j=1}^{i-1}L_j:1:L_i-1+\sum _{j=1}^{i-1}L_j$ $($if $i=1$, $\sum _{j=1}^{i-1}L_j=0$$)$,
                    solve the following linear systems on each time level
                        \begin{numcases}{}
                            (\frac{\alpha_{\nu}}{\triangle t}I+\frac{\beta_{\nu}}{\omega}A)x_{n+\nu}^{(k)}=
                            \nonumber\\\qquad
                            \sum \limits_{j=0}^{\nu}\beta_{j}((\frac{1}{\omega}-1)Ax_{n+j}^{(k-1)}-By_{n+j}^{(k-1)}+f_{n+j})
                            -\sum \limits_{j=0}^{\nu-1}
                            (\frac{\alpha_{j}}{\triangle t}I+\frac{\beta_{j}}{\omega}A)x_{n+j}^{(k)},\nonumber\\
                            \frac{\beta_{\nu}}{\tau}Qy_{n+\nu}^{(k)} = \sum \limits_{j=0}^{\nu}\beta_{j}
                            (B^Tx_{n+j}^{(k)}+\frac{1}{\tau}Qy_{n+j}^{(k-1)}+g_{n+j}) -
                            \sum \limits_{j=0}^{\nu-1} \frac{\beta_{j}}{\tau}Qy_{n+j}^{(k)}.\nonumber
                        \end{numcases}
                    \item[End]
                \end{itemize}
            \item[End]
        \end{itemize}
    \item[End]
\end{itemize}
\end{method}

\section{Numerical Results}
\label{sec-num}

In this section, numerical tests are performed to demonstrate the correctness
of theoretical results presented in previous sections and the efficiency of the
DABSOR method with windowing technique, i.e. the Method \ref{dwr-gsor-win},
for solving the linear DAEs derived from time-dependent Stokes equations.

Consider the two-dimensional time-dependent
Stokes equations on the domain $\Omega=\{-1\le
{\mbox{x}}\le 1, -1\le {\mbox{y}}\le 1\}$:
\begin{eqnarray}
\label{time-stokes-2d}
\left\{
    \begin{array}{l}
    \frac{\partial u}{\partial t}-
    \nu(\frac{\partial^2 u}{\partial{\mbox{x}}^2}+
    \frac{\partial^2 u }{\partial {\mbox{y}}^2})+
    \frac{\partial p}{\partial {\mbox{x}}}=0,\\
    \frac{\partial v}{\partial t}-
    \nu(\frac{\partial^2 v}{\partial{\mbox{x}}^2}+
    \frac{\partial^2 v}{\partial {\mbox{y}}^2})+
    \frac{\partial p}{\partial {\mbox{y}}}=0,\\
    \frac{\partial u}{\partial {\mbox{x}}}+
    \frac{\partial v}
    {\partial {\mbox{y}}}=0.
    \end{array}
    \right.
\end{eqnarray}
Note that the analytic solution of equations (\ref{time-stokes-2d}) is
of the form
\begin{eqnarray*}
\left\{
    \begin{array}{l}
    u_{\star}(\mbox{x},\mbox{y};t)=\mbox{\bf u}
    ({\mbox{y}})e^{\theta{\mbox{x}}-\zeta t},\quad
    \mbox{\bf u}({\mbox{y}})=
    c_1\sin(\theta{\mbox{y}})+
    \frac{2\kappa}{\theta}c_2
    \sin(\kappa {\mbox{y}}),\\
    v_{\star}(\mbox{x},\mbox{y};t)=\mbox{\bf v}
    ({\mbox{y}})e^{\theta{\mbox{x}}-\zeta t},\quad
    \mbox{\bf v}({\mbox{y}})=
    c_1\cos(\theta{\mbox{y}})+
    2 c_2\cos(\kappa {\mbox{y}}),\\
    p_{\star}(\mbox{x},\mbox{y};t)=\mbox{\bf p}
    ({\mbox{y}})e^{\theta{\mbox{x}}-\zeta t},\quad
    \mbox{\bf p}({\mbox{y}})=
    \frac{\zeta}{\theta}c_1\sin(\theta{\mbox{y}}),
    \end{array}
    \right.
\end{eqnarray*}
where the parameters are chosen the same as in \cite{BaiYang12}, i.e.
set $\nu=1$, $\theta=1$, $\zeta=11.6348$, $\kappa=3.5545$,
$c_1=3.390472650419484$, $c_2=1$.

Since this paper focuses on the study of solving linear DAEs by waveform relaxation
methods, less attention is paid to spacial discretization of the
time-dependent Stokes equations (\ref{time-stokes-2d}). For a uniform spacial
grid with stepsize $h_{\mbox{x}}=\frac{2}{\ell_{\mbox{\tiny x}}+1}$ and
$h_{\mbox{y}}=\frac{2}{\ell_{\mbox{\tiny y}}+1}$, we simply choose Scheme II
defined in \cite{BaiYang12} for the implementation of the DABSOR method,
which applies the centered difference scheme to the Laplacian, performs the forward
difference scheme to the pressure variable, and discretizes the third equation
in (\ref{time-stokes-2d}) by the backward difference scheme. Then we obtain a linear
DAEs of the form (\ref{dis-dae-LCC}), the details of its coefficient matrices
and right-hand side vector-valued function can be found in \cite{BaiYang12}.
Besides, the choices of the preconditioner matrix $Q$ are shown in Table \ref{tab:caseQ}.
\begin{table}[hbp]
\renewcommand{\arraystretch}{1.5}
\setlength{\abovecaptionskip}{0pt}
  \setlength{\belowcaptionskip}{10pt}
  \centering{
\caption{\label{tab:caseQ}
The Choices of the Preconditioner Matrix $Q$}
\begin{tabular}{|c|c|c|}
 \hline
 Case No. & Matrix $Q$ & Description \\
\hline
 $Q_1$ & $B^T \widehat A^{-1} B $ &
 $\widehat A = \mbox{tridiag}(A)$\\
\hline
 $Q_2$ &$B^T \widehat A^{-1} B $ &
 $\widehat A = \mbox{diag}(A)$ \\
\hline
\end{tabular}}
\end{table}

Due to the stiffness of the linear DAEs, the backward differentiation formulae
(BDF) of order 1 to order 6 are selected to be the linear multistep
formulae for the DABSOR method. The coefficients of backward differentiation
formulae are shown in Table \ref{tab:BDF-coe}.
\begin{table}
\renewcommand{\arraystretch}{1.4}
\setlength{\abovecaptionskip}{0pt}
  \setlength{\belowcaptionskip}{10pt}
  \centering{
\caption{\label{tab:BDF-coe}
The coefficients of BDF}
\begin{tabular}{|c||c|c|c|c|c|c|c|c|}
 \hline
 Order $\nu$ & $\beta_{\nu}$ &
 $\alpha_6$ & $\alpha_5$ & $\alpha_4$ & $\alpha_3$ & $\alpha_2$ & $\alpha_1$ & $\alpha_0$ \\\hline\hline
 $1$ & $1$ &
 &  &  &  &  & $1$ & $-1$ \\\hline
 $2$ & $\frac{2}{3}$ &
 &  &  &  & $1$ & $-\frac{4}{3}$ & $\frac{1}{3}$ \\\hline
 $3$ & $\frac{6}{11}$ &
 &  &  & $1$ & $-\frac{8}{11}$ & $\frac{9}{11}$ & $-\frac{2}{11}$ \\\hline
 $4$ & $\frac{12}{25}$ &
 &  & $1$ & $-\frac{48}{25}$ & $\frac{36}{25}$ & $-\frac{16}{25}$ & $\frac{3}{25}$ \\\hline
 $5$ & $\frac{60}{137}$ &
 & $1$ & $-\frac{300}{137}$ & $\frac{300}{137}$ & $-\frac{200}{137}$ & $\frac{75}{137}$ & $-\frac{12}{137}$ \\\hline
 $6$ & $\frac{60}{147}$ &
 $1$ & $-\frac{360}{147}$ & $\frac{450}{147}$ & $-\frac{400}{147}$ & $\frac{225}{147}$ & $-\frac{72}{147}$ & $\frac{10}{147}$ \\\hline
\end{tabular}}
\end{table}

In fact, the purpose for simulation is to
compute the approximate solution of the time-dependent Stokes
equations (\ref{time-stokes-2d}) on a finite time interval
$\Omega_t=\cup _{i=1}^N \Omega_t^{(i)}$, where
$\Omega_t^{(i)}=(\txt{T}_1+\triangle t \times \sum _{j=1}^{i-1} L_j,
\txt{T}_1+\triangle t \times \sum _{j=1}^{i} L_j]$ is the $i$-th
window of the DABSOR method (if $i=1$, $\sum _{j=1}^{i-1} L_j=0$).
Here, $\triangle t$ represents the time stepsize, $N$ denotes the number
of windows, and $L_i$ is the number of time steps on the $i$-th window.
The stopping criterion on the $i$-th window of the DABSOR method
is set to be
\begin{eqnarray}
\label{direct-residual}
\epsilon ^{(k,i)} =
\frac{
\sup _{\Omega\times\Omega_t^{(i)}}
\left\{\bigl|u_h^{(k,i)}-u_h^{(*,i)}\bigr|,\
\bigl|v_h^{(k,i)}-v_h^{(*,i)}\bigr|,\
\bigl|p_h^{(k,i)}-p_h^{(*,i)}\bigr| \right\}
}{
\sup _{\Omega\times\Omega_t^{(i)}}
\left\{\bigl|u_h^{(*,i)}\bigr|,\
\bigl|v_h^{(*,i)}\bigr|,\
\bigl|p_h^{(*,i)}\bigr| \right\}
} < 10^{-6},
\end{eqnarray}
where $u_h^{(k,i)}$, $v_h^{(k,i)}$, $p_h^{(k,i)}$ are the $k$-th
iterate on the $i$-th window of the DABSOR method,
and $u_h^{(*,i)}$, $v_h^{(*,i)}$, $p_h^{(*,i)}$ are the entries
of the exact solution on the $i$-th window of the linear
DAEs (\ref{dis-dae-LCC}) derived from equations (\ref{time-stokes-2d}).
Moreover, the initial waves are chosen to be
\begin{align*}
&u_h^{(0)}=\frac{1}{1+10000\zeta t}\,
u_{\star}(\mbox{x},\mbox{y};0),\\
&v_h^{(0)}=\frac{1}{1+10000\zeta t}\,
v_{\star}(\mbox{x},\mbox{y};0)
\intertext{and}
&p_h^{(0)}=\frac{1}{1+10000\zeta t}\,
p_{\star}(\mbox{x},\mbox{y};0).
\end{align*}

Since the DABSOR method is integrated with windowing
technique, long time interval simulation of the time-dependent Stokes equations
(\ref{time-stokes-2d}) can be obtained by simply
adding as many windows as required to the end of the existing time
interval. Hence, it is not necessary to choose a long time
interval for numerical tests. In the sequel, the time step
is fixed as $\triangle t=0.001$,
and the simulation time interval of the time-dependent Stokes equations
(\ref{time-stokes-2d}) is $(0.01,0.13]$. Due to
the application of linear multistep formulae to the DABSOR method,
the exact solution of the time-dependent Stokes equations (\ref{time-stokes-2d})
on $[0,0.01]$ is taken to serve as the initial values.
For a precise and comprehensive demonstration, the following three
subsections present the numerical results in three different aspects.

\subsection{Optimal Convergence Curve}
\label{sec-num-occ}
In Theorem \ref{th-dis-opt} and Remark \ref{rk-opt-con-curve}, there is a curve
in $\omega$-$\tau$ plane related to the optimality of the CABSOR method called
the optimal convergence curve. The curve is shown in this section in
different situations.

\begin{figure}
            \centering
            \begin{tabular}{cc}
                \includegraphics[scale=0.4]{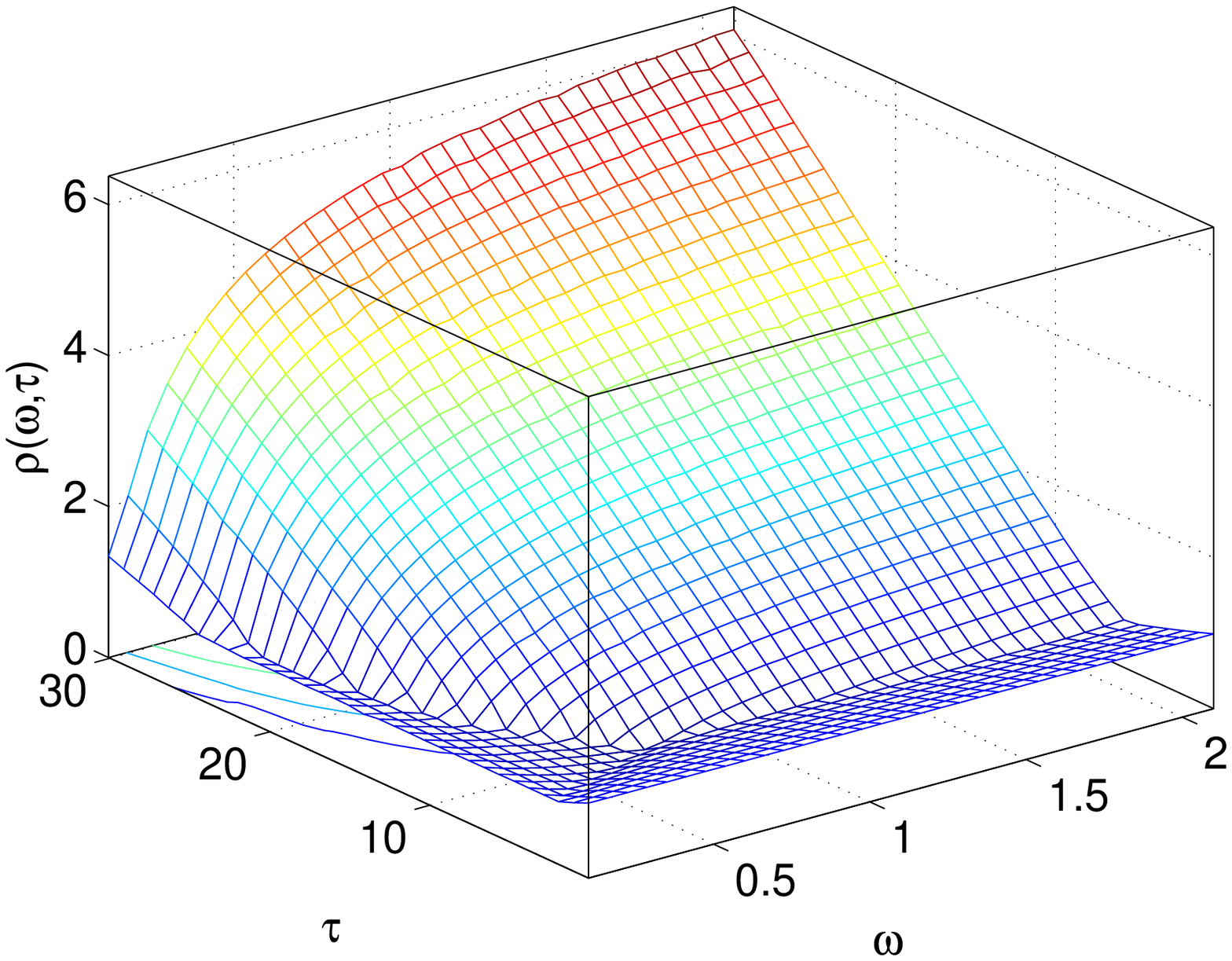}
                &
                \includegraphics[scale=0.4]{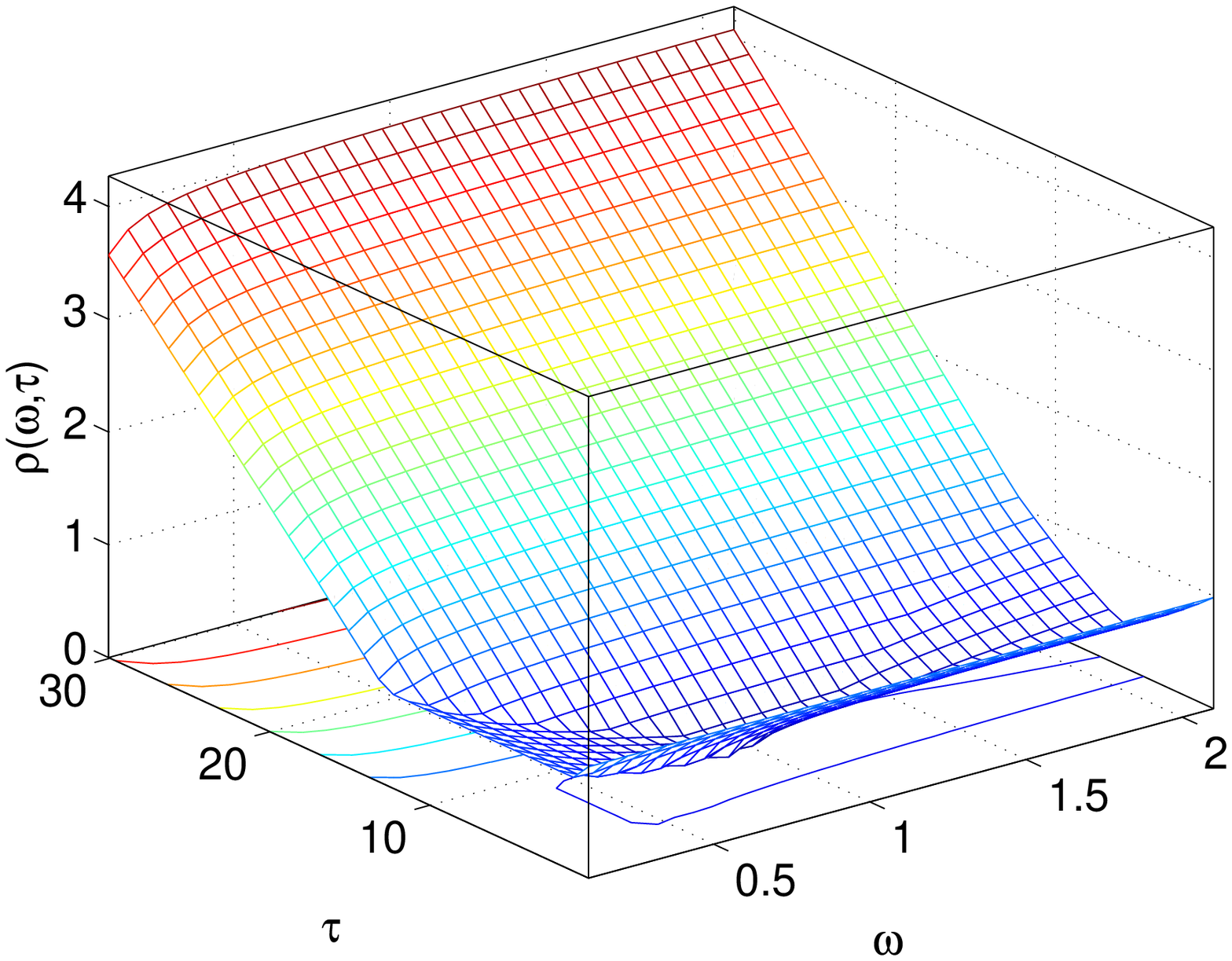}
                \\
                (a) & (b)\\
            \end{tabular}
            \caption{Surfaces of spectral radii of linear convolution operator
            $\mathcal{K}_{\text{\tiny$\triangle t$}}$ on finite time interval
            with respect to $\omega$ and $\tau$ based on BDF(1), $12\times 12$ grid
            and preconditioners (a) $Q_1$ and (b) $Q_2$.
            }
            \label{fig:dwr-radii-order1}
\end{figure}

\begin{figure}
            \centering
            \begin{tabular}{cc}
                \includegraphics[scale=0.4]{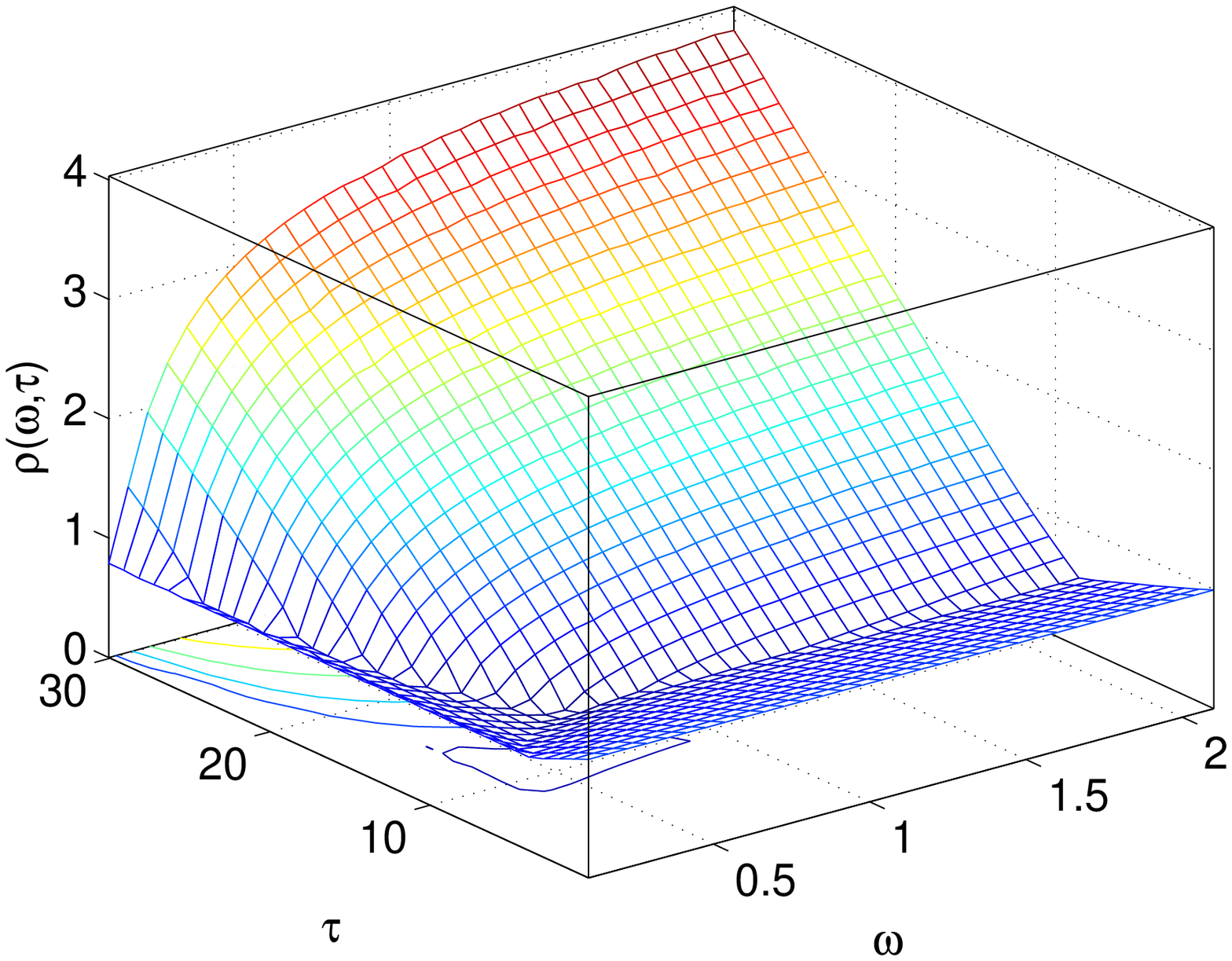}
                &
                \includegraphics[scale=0.4]{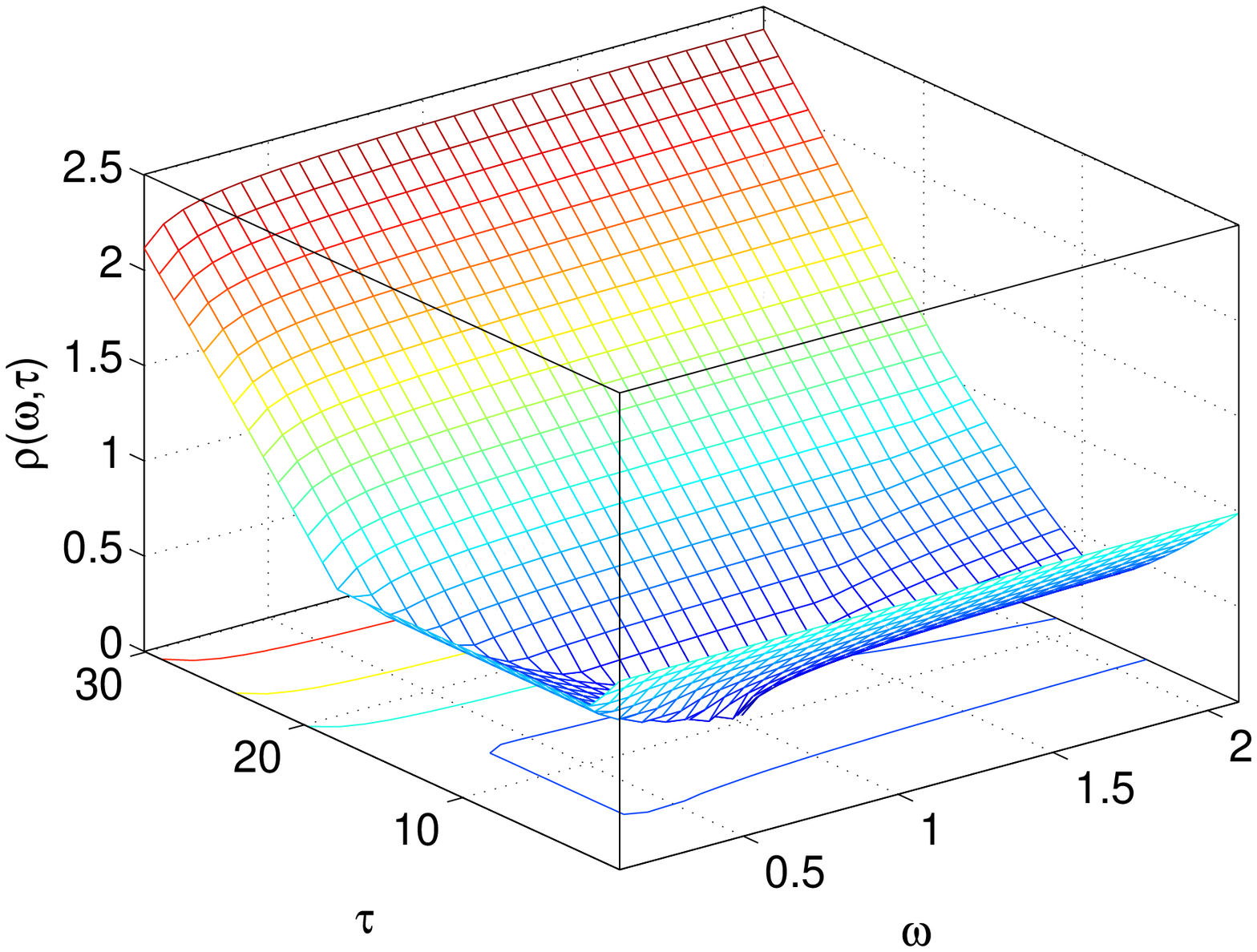}
                \\
                (a) & (b)\\
            \end{tabular}
            \caption{Surfaces of spectral radii of linear convolution operator
            $\mathcal{K}_{\text{\tiny$\triangle t$}}$ on finite time interval
            with respect to $\omega$ and $\tau$ based on BDF(2), $12\times 12$ grid
            and preconditioners (a) $Q_1$ and (b) $Q_2$.
            }
            \label{fig:dwr-radii-order2}
\end{figure}

\begin{figure}
            \centering
            \begin{tabular}{cc}
                \includegraphics[scale=0.4]{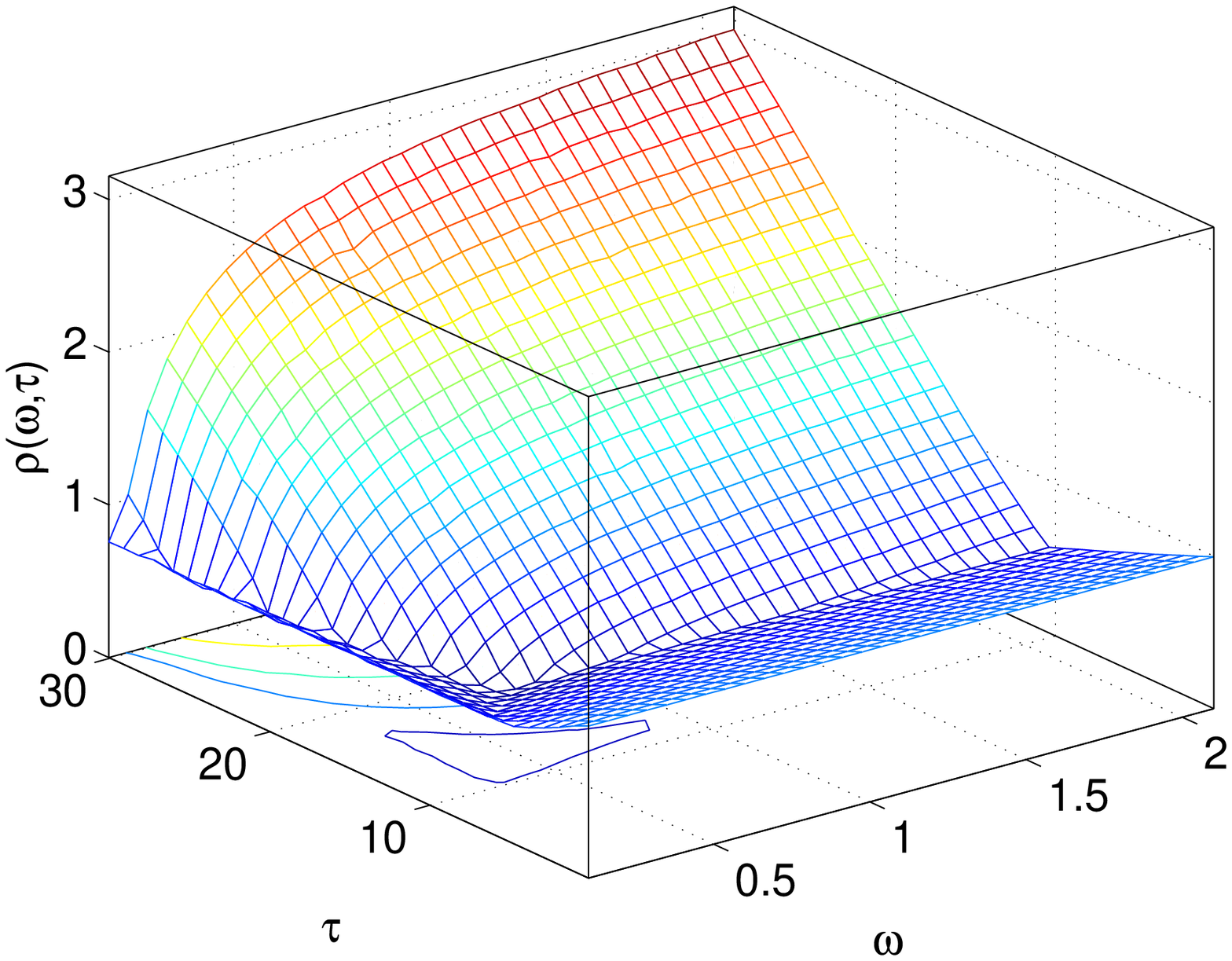}
                &
                \includegraphics[scale=0.4]{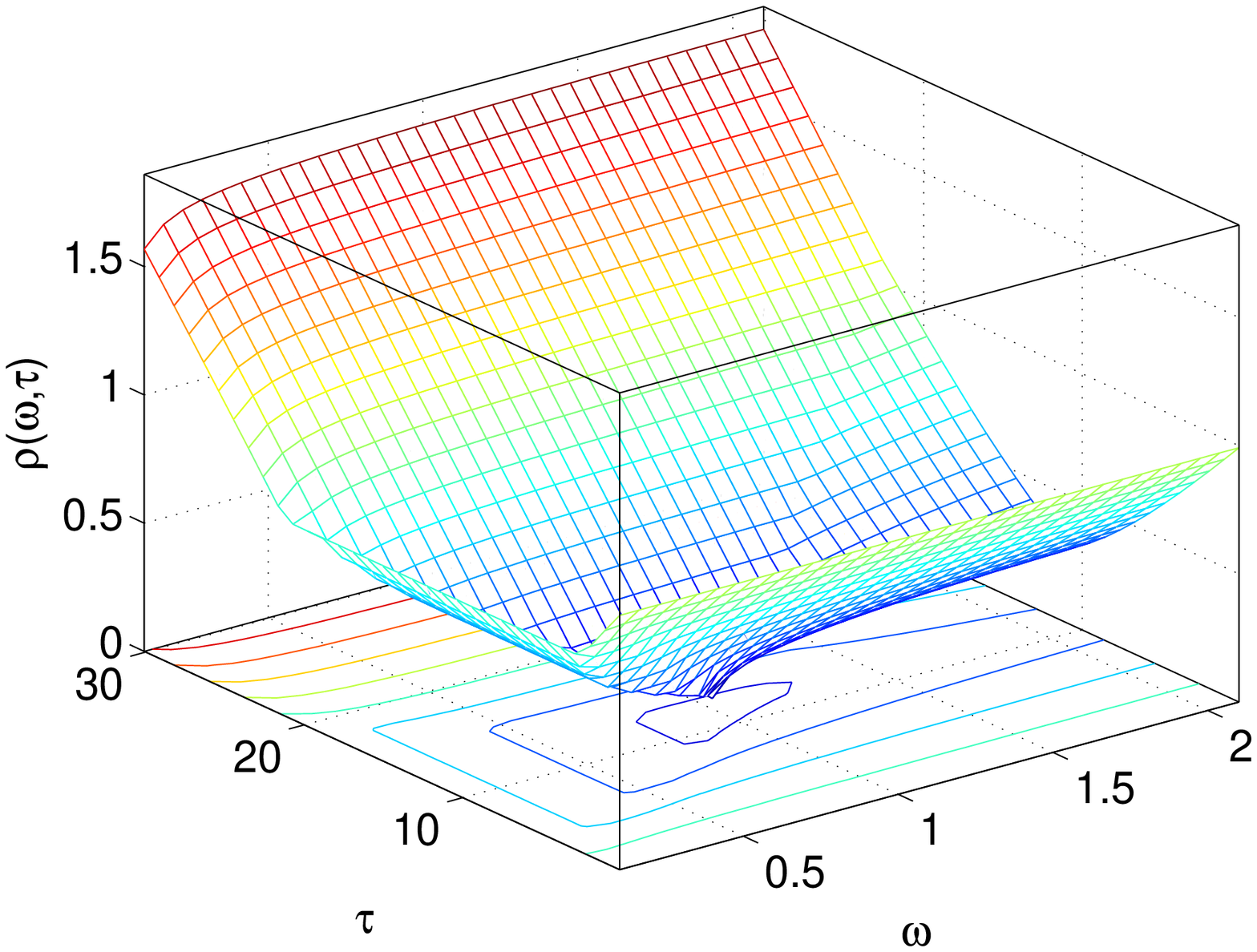}
                \\
                (a) & (b)\\
            \end{tabular}
            \caption{Surfaces of spectral radii of linear convolution operator
            $\mathcal{K}_{\text{\tiny$\triangle t$}}$ on finite time interval
            with respect to $\omega$ and $\tau$ based on BDF(3), $12\times 12$ grid
            and preconditioners (a) $Q_1$ and (b) $Q_2$.
            }
            \label{fig:dwr-radii-order3}
\end{figure}

\begin{figure}
            \centering
            \begin{tabular}{cc}
                \includegraphics[scale=0.4]{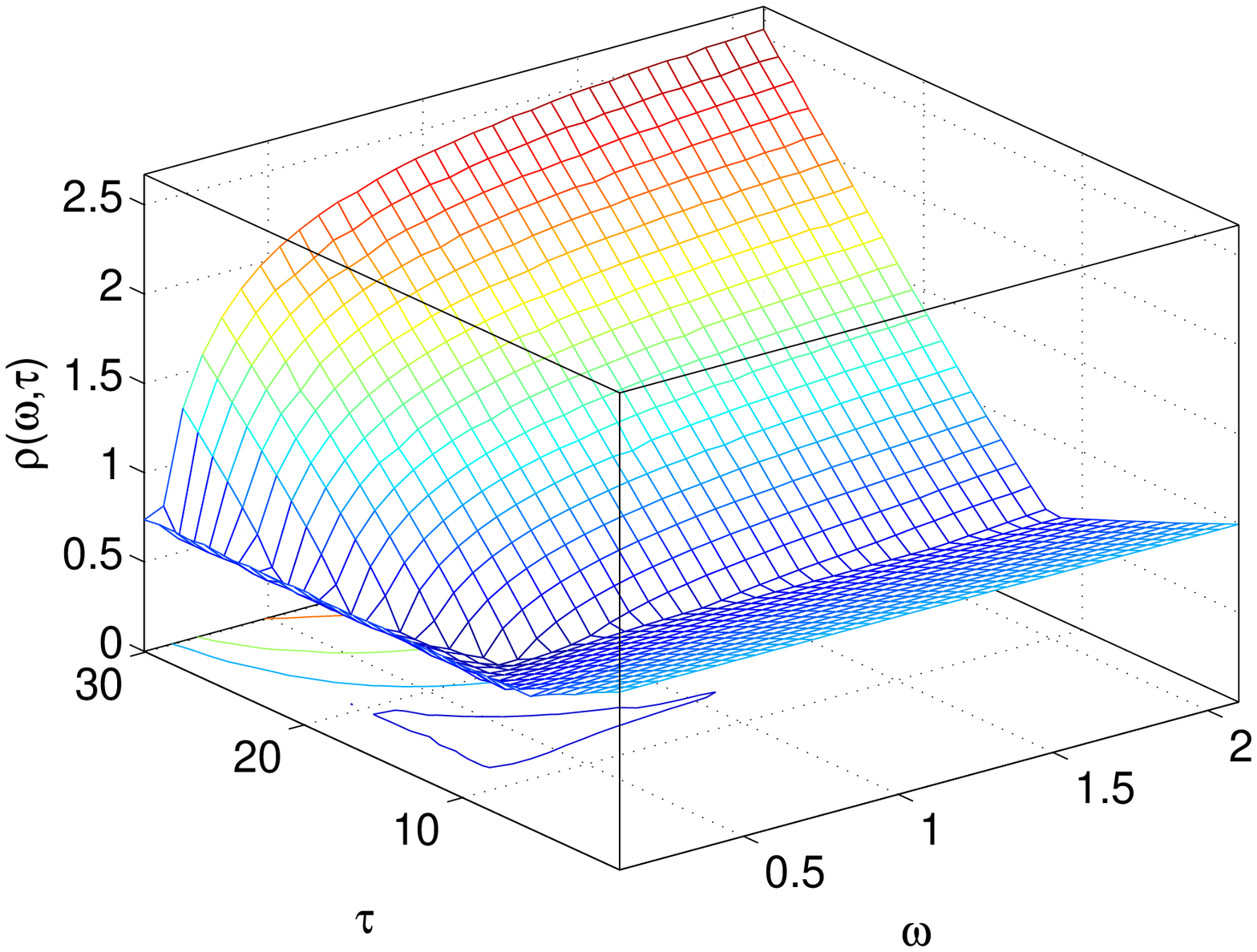}
                &
                \includegraphics[scale=0.4]{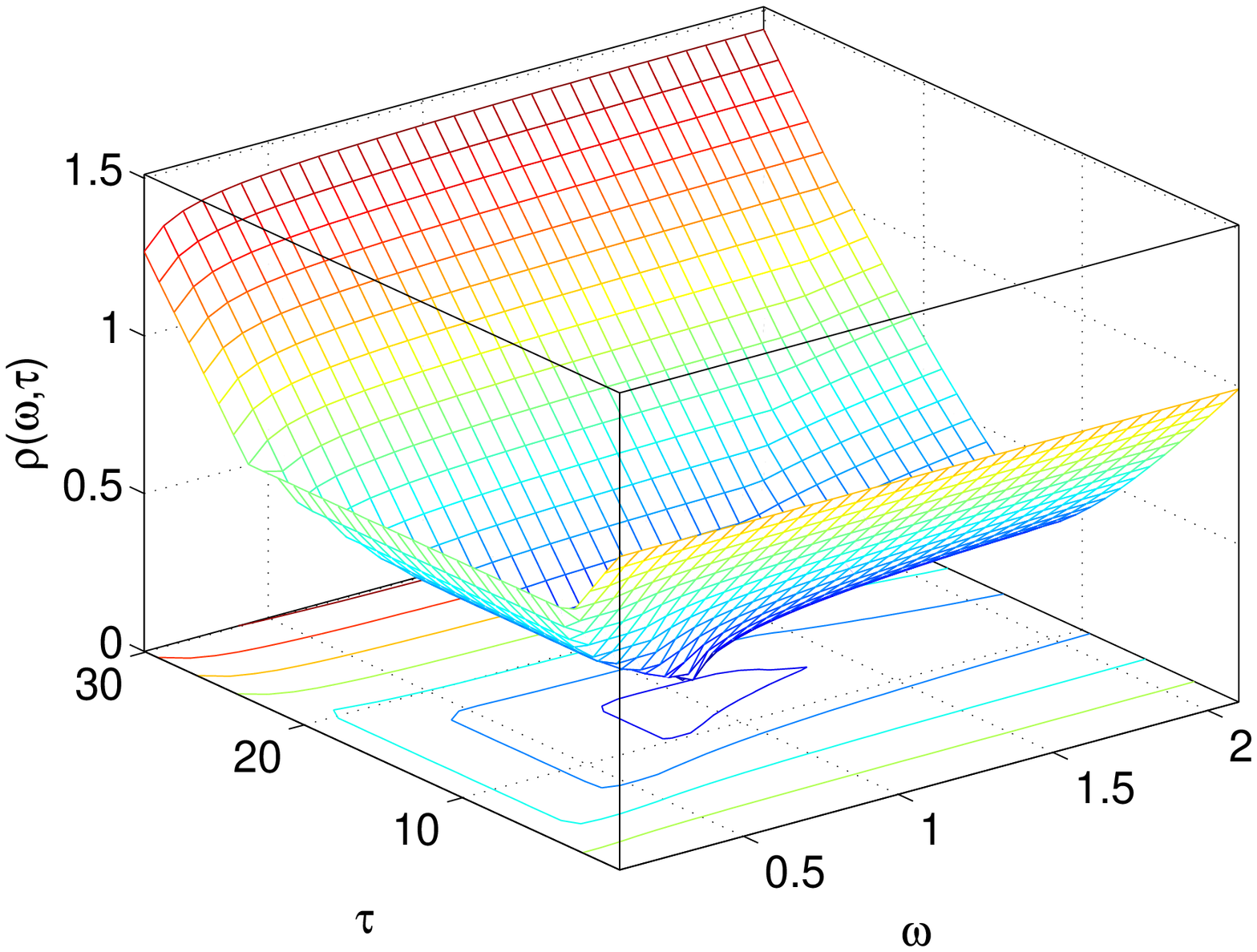}
                \\
                (a) & (b)\\
            \end{tabular}
            \caption{Surfaces of spectral radii of linear convolution operator
            $\mathcal{K}_{\text{\tiny$\triangle t$}}$ on finite time interval
            with respect to $\omega$ and $\tau$ based on BDF(4), $12\times 12$ grid
            and preconditioners (a) $Q_1$ and (b) $Q_2$.
            }
            \label{fig:dwr-radii-order4}
\end{figure}

\begin{figure}
            \centering
            \begin{tabular}{cc}
                \includegraphics[scale=0.4]{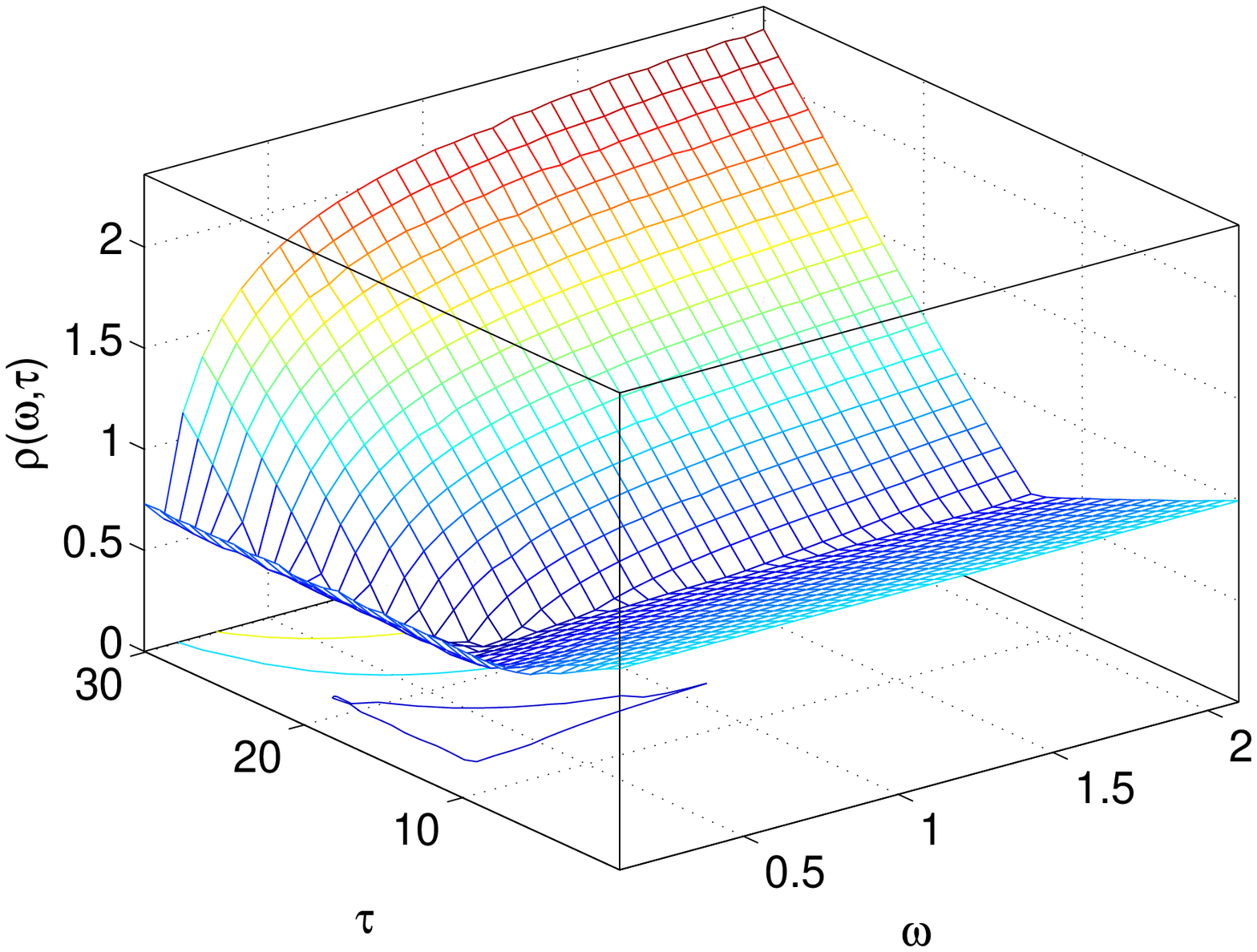}
                &
                \includegraphics[scale=0.4]{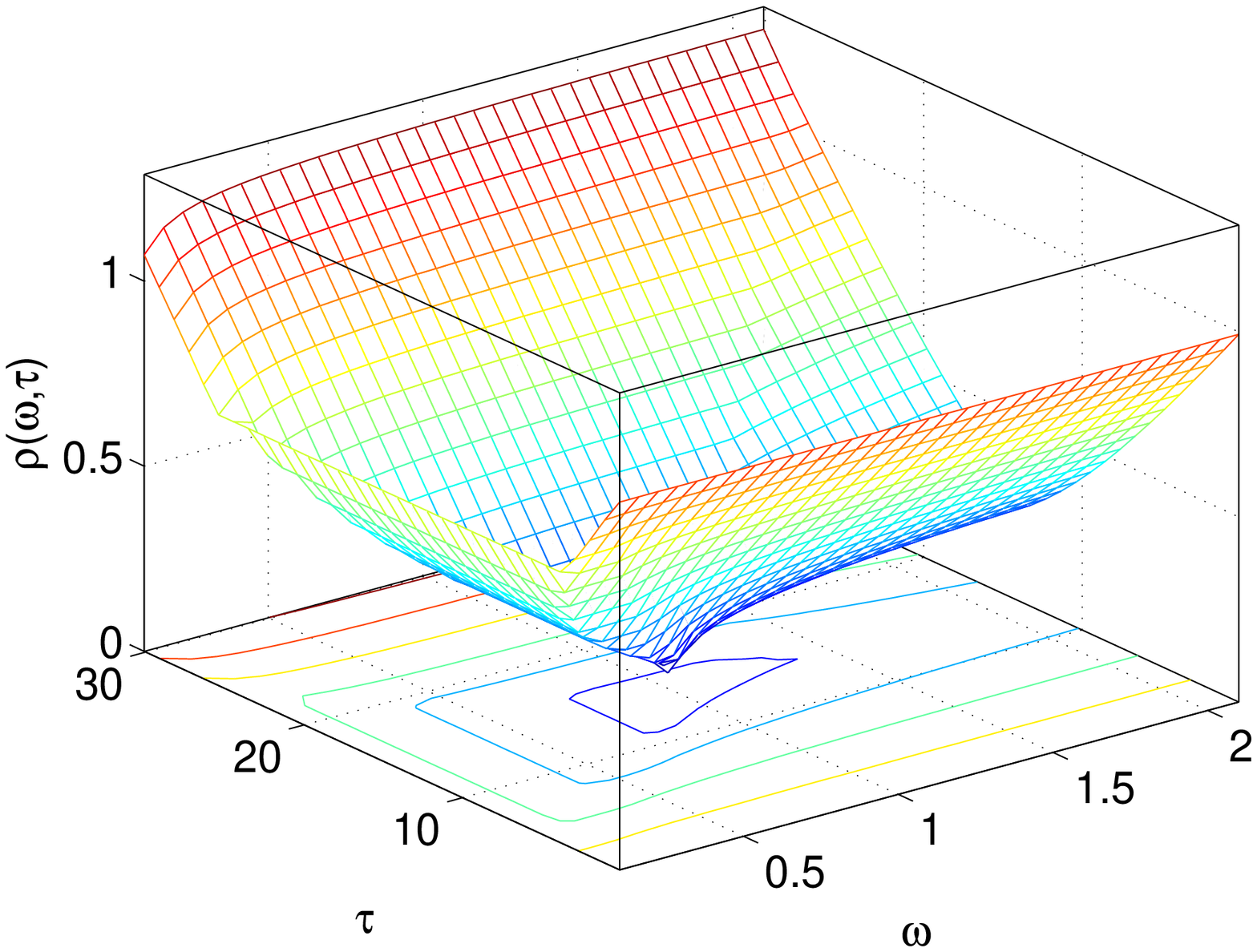}
                \\
                (a) & (b)\\
            \end{tabular}
            \caption{Surfaces of spectral radii of linear convolution operator
            $\mathcal{K}_{\text{\tiny$\triangle t$}}$ on finite time interval
            with respect to $\omega$ and $\tau$ based on BDF(5), $12\times 12$ grid
            and preconditioners (a) $Q_1$ and (b) $Q_2$.
            }
            \label{fig:dwr-radii-order5}
\end{figure}

\begin{figure}
            \centering
            \begin{tabular}{cc}
                \includegraphics[scale=0.4]{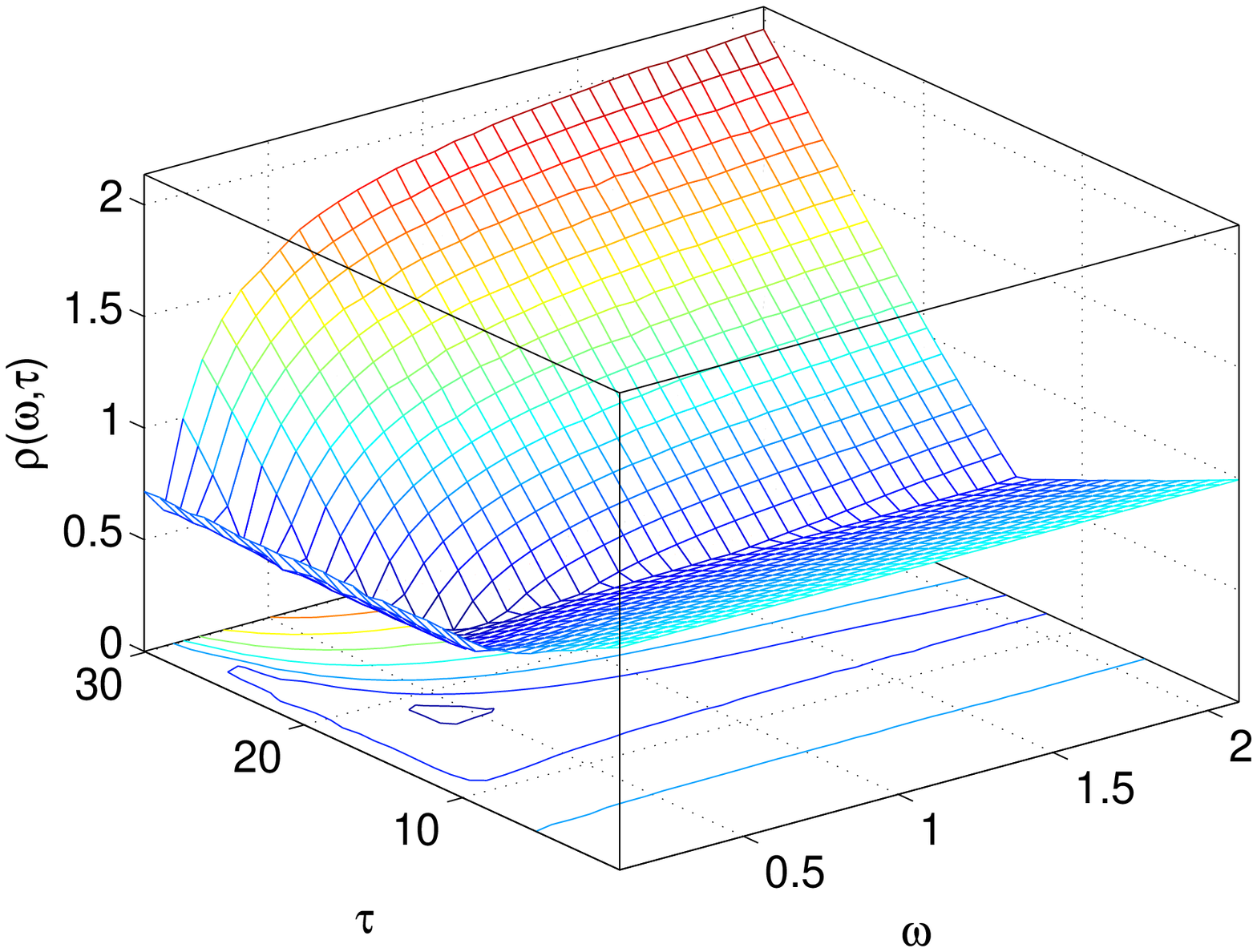}
                &
                \includegraphics[scale=0.4]{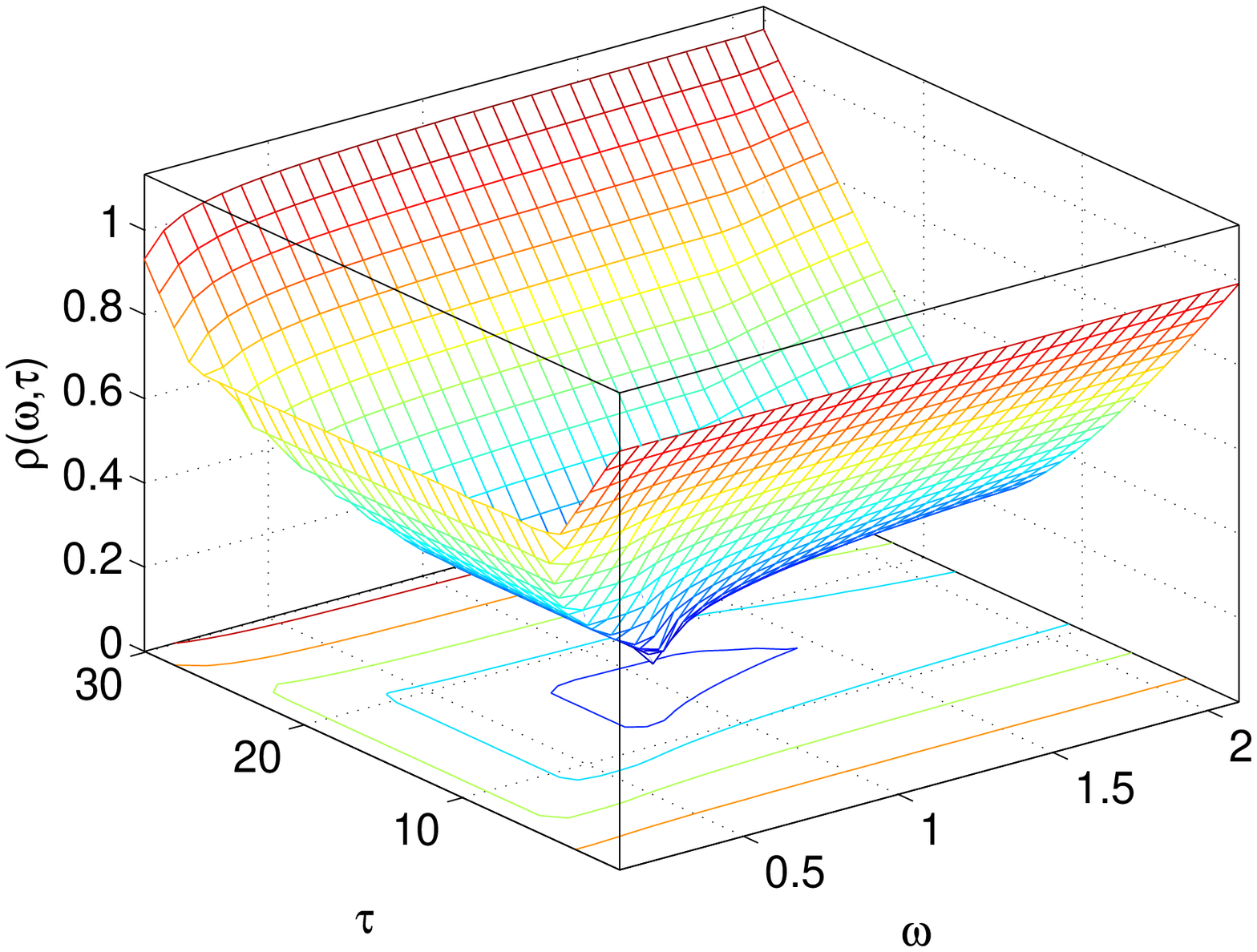}
                \\
                (a) & (b)\\
            \end{tabular}
            \caption{Surfaces of spectral radii of linear convolution operator
            $\mathcal{K}_{\text{\tiny$\triangle t$}}$ on finite time interval
            with respect to $\omega$ and $\tau$ based on BDF(6), $12\times 12$ grid
            and preconditioners (a) $Q_1$ and (b) $Q_2$.
            }
            \label{fig:dwr-radii-order6}
\end{figure}

Surfaces of spectral radii of the discrete linear convolution operator
$\mathcal{K}_{\text{\tiny$\triangle t$}}$ based on six different linear multistep
formulae like BDF(1-6), one grid size as $12\times 12$ and two different choices of
preconditioners $Q$ in Table \ref{tab:caseQ} are shown in
Figures \ref{fig:dwr-radii-order1}-\ref{fig:dwr-radii-order6}. Here, BDF($i$)
represents the backward differentiation formula of order $i$. According to
Theorem \ref{th-dis-opt} and Remark \ref{rk-opt-con-curve}, the optimal iterative
parameter pair $(\omega_{\opt},\tau _{\opt})$ is not a single point in
$\omega$-$\tau$ plane, all possible choices of optimal iterative parameter
pair $(\omega_{\opt},\tau_{\opt})$ lead to a finite length parameterized curve
with respect to $\delta\in(\delta_{\min},\delta_{\max})$, i.e. the optimal
convergence curve.

After careful observation of the surfaces based on
preconditioner $Q_1$ in Figures \ref{fig:dwr-radii-order1}-\ref{fig:dwr-radii-order6},
we find that the lower bound of spectral radii of the discrete linear convolution
operator $\mathcal{K}_{\text{\tiny$\triangle t$}}$ shown in each surface based
on BDF of six different orders is available along a 3-D curve of certain length.
Obviously, the projection of such 3-D curve to the $\omega$-$\tau$ plane is just
the optimal convergence curve, which coincides with the description in
Theorem \ref{th-dis-opt} and Remark \ref{rk-opt-con-curve}. However, for the case
of preconditioner $Q_2$, the length of optimal convergence curve decreases sharply
when the order of BDF increases. Especially for BDF(4-6), the optimal
convergence curve shrinks to a single point. It means that the optimal
convergence factor of the DABSOR method based on preconditioner $Q_2$ is
much more sensitive to the choice of iterative parameters than that based on preconditioner $Q_1$.
According to the expressions
of optimal iterative parameters $\omega_{\opt}$ and $\tau_{\opt}$, we find that
the length of optimal convergence curve closely related to matrices $A$, $B$ and
preconditioner $Q$. Apparently, the difference between surfaces in each figure
is caused by choosing different preconditioner $Q$. There are also good
news for preconditioner $Q_2$, that is, the lower bound of spectral radii of
the discrete linear convolution operator $\mathcal{K}_{\text{\tiny$\triangle t$}}$ based on
preconditioner $Q_2$ is much smaller than that based on preconditioner $Q_1$. It means
that the DABSOR method with optimal convergence parameters based on preconditioner $Q_2$
is much faster than that based on preconditioner $Q_1$.
Therefore, the preconditioner matrix $Q$ should be chosen carefully.
\subsection{Optimal Convergence Factor: Theoretical vs Practical}
\label{sec-num-ocf}
In Theorems \ref{th-dis-con-fin} and \ref{th-dis-opt}, the general spectral
radius formula of discrete-time waveform relaxation methods for solving general
linear DAEs and the optimal convergence factor formula of the DABSOR method
for solving linear DAEs derived from the time-dependent Stokes equations are
discovered on finite time interval respectively. In this subsection,
the comparison between theoretical and practical value of optimal convergence
factor of the DABSOR method is exhibited. The computation of
theoretical value is based on Theorems \ref{th-dis-con-fin} and
\ref{th-dis-opt}. The practical value represents the average experimental
convergence rate.

The comparisons between theoretical and practical value of optimal convergence
factor based on six different linear multistep formulae like BDF(1-6),
two grid sizes as $12\times 12$, $24\times 24$ and two different choices of
preconditioners $Q$ in Table \ref{tab:caseQ} are presented in Tables
\ref{tab:th_mea_size12_Q1}-\ref{tab:th_mea_size24_Q2}. In these tables, DTOCF
denotes the theoretical value of optimal convergence factor of the DABSOR method
on finite time interval, APOCF represents the practical value of static iterative
method for solving the system of linear equations with respect to coefficient
matrix $\mathcal{A}$ in linear DAEs (\ref{dis-dae-LCC}) and iterative matrix
$\bbb{K}(\sigma)$ with $\omega_{\opt}$ and $\tau_{\opt}$, and DPOCF($N$)
is the practical value of optimal convergence factor of the DABSOR method
with $N$ windows. After observing Tables
\ref{tab:th_mea_size12_Q1}-\ref{tab:th_mea_size24_Q2}, we find that
the practical value DPOCF($N$) of the DABSOR method decreases when the
the number of windows increases, and the larger the number of windows,
the closer the practical value DPOCF($N$) to the theoretical value DTOCF.
Specifically, for smaller number of windows, the computation load of
the DABSOR method increases intensively because of larger spurious oscillations
discussed in section \ref{subsec-dis-dwr-wintec} occurred on each window. Thus,
the practical value DPOCF($N$) with small number of windows is far beyond
the theoretical value DTOCF. When $N$ the number of windows increases,
the spurious oscillations on each window becomes less apparent, and
the DABSOR method tends to be more efficient. In addition, the practical
value APOCF of the static iterative method is always smaller than the
theoretical value DTOCF and practical value DPOCF($N$) of the DABSOR method,
indicating that the convergence rate of the DABSOR method is unlikely to
be faster than the corresponding static iterative method. On the other hand,
theoretical and practical values based on preconditioner $Q_2$ are much
smaller than that based on preconditioner $Q_1$, which means that
preconditioner $Q_2$ always leads to faster convergence rate. The findings
are consist with the results in subsection \ref{sec-num-occ}.

\begin{table}
\setlength{\abovecaptionskip}{0pt}
\setlength{\belowcaptionskip}{10pt}
\centering{
\caption{\label{tab:th_mea_size12_Q1} Theoretical vs Practical: the optimal convergence factor
based on $12\times 12$ grid and preconditioner $Q_1$}
\begin{tabular}{|l|c|c|c|c|c|c|}\hline
 & BDF(1) & BDF(2) & BDF(3) & BDF(4) & BDF(5) & BDF(6) \\ \hline
 DTOCF     & 0.3590 & 0.3921 & 0.4991 & 0.4812 & 0.5058 & 0.4943 \\ \hline
 APOCF     & 0.2154 & 0.2512 & 0.2154 & 0.2154 & 0.2154 & 0.2154 \\ \hline
 DPOCF(6)  & 0.7252 & 0.7834 & 0.8555 & 0.8918 & 0.9041 & 0.9056 \\ \hline
 DPOCF(12) & 0.6519 & 0.6679 & 0.7960 & 0.8287 & 0.8417 & 0.8414 \\ \hline
 DPOCF(20) & 0.6138 & 0.5695 & 0.7321 & 0.7687 & 0.7832 & 0.7803 \\ \hline
 DPOCF(30) & 0.4527 & 0.4808 & 0.5985 & 0.5950 & 0.6152 & 0.6065 \\ \hline
 DPOCF(40) & 0.3915 & 0.4206 & 0.5501 & 0.5424 & 0.5649 & 0.5530 \\ \hline
 DPOCF(60) & 0.3237 & 0.3435 & 0.4855 & 0.4737 & 0.4946 & 0.4803 \\ \hline
\end{tabular}}
\end{table}
\begin{table}
\setlength{\abovecaptionskip}{0pt}
\setlength{\belowcaptionskip}{10pt}
\centering{
\caption{\label{tab:th_mea_size12_Q2} Theoretical vs Practical: the optimal convergence factor
based on $12\times 12$ grid and preconditioner $Q_2$}
\begin{tabular}{|l|c|c|c|c|c|c|}\hline
 & BDF(1) & BDF(2) & BDF(3) & BDF(4) & BDF(5) & BDF(6) \\ \hline
 DTOCF     & 0.2940 & 0.1804 & 0.0559 & 0.1098 & 0.1250 & 0.1084 \\ \hline
 APOCF     & 0.0631 & 0.0316 & 0.0316 & 0.0100 & 0.0316 & 0.0100 \\ \hline
 DPOCF(6)  & 0.6081 & 0.6302 & 0.6661 & 0.6585 & 0.7341 & 0.7282 \\ \hline
 DPOCF(12) & 0.4892 & 0.4762 & 0.4728 & 0.4603 & 0.5754 & 0.5601 \\ \hline
 DPOCF(20) & 0.4115 & 0.3490 & 0.3300 & 0.3051 & 0.4420 & 0.4215 \\ \hline
 DPOCF(30) & 0.3516 & 0.2673 & 0.2053 & 0.2107 & 0.2457 & 0.2205 \\ \hline
 DPOCF(40) & 0.3006 & 0.2083 & 0.1493 & 0.1590 & 0.1838 & 0.1520 \\ \hline
 DPOCF(60) & 0.2379 & 0.1412 & 0.0926 & 0.0895 & 0.1176 & 0.0827 \\ \hline
\end{tabular}}
\end{table}
\begin{table}
\setlength{\abovecaptionskip}{0pt}
\setlength{\belowcaptionskip}{10pt}
\centering{
\caption{\label{tab:th_mea_size24_Q1} Theoretical vs Practical: the optimal convergence factor
based on $24\times 24$ grid and preconditioner $Q_1$}
\begin{tabular}{|l|c|c|c|c|c|c|}\hline
 & BDF(1) & BDF(2) & BDF(3) & BDF(4) & BDF(5) & BDF(6) \\ \hline
 DTOCF     & 0.3632 & 0.3298 & 0.5236 & 0.3408 & 0.5004 & 0.3493 \\ \hline
 APOCF     & 0.2154 & 0.2154 & 0.2154 & 0.2154 & 0.2154 & 0.2154 \\ \hline
 DPOCF(6)  & 0.6811 & 0.7525 & 0.8203 & 0.8188 & 0.8808 & 0.8567 \\ \hline
 DPOCF(12) & 0.5875 & 0.6930 & 0.7603 & 0.7123 & 0.8061 & 0.7616 \\ \hline
 DPOCF(20) & 0.5860 & 0.6271 & 0.6926 & 0.6156 & 0.7365 & 0.6722 \\ \hline
 DPOCF(30) & 0.4684 & 0.4844 & 0.6039 & 0.4861 & 0.6026 & 0.4848 \\ \hline
 DPOCF(40) & 0.4065 & 0.4189 & 0.5589 & 0.4168 & 0.5514 & 0.4158 \\ \hline
 DPOCF(60) & 0.3412 & 0.3369 & 0.4973 & 0.3323 & 0.4816 & 0.3284 \\ \hline
\end{tabular}}
\end{table}
\begin{table}
\setlength{\abovecaptionskip}{0pt}
\setlength{\belowcaptionskip}{10pt}
\centering{
\caption{\label{tab:th_mea_size24_Q2} Theoretical vs Practical: the optimal convergence factor
based on $24 \times 24$ grid and preconditioner $Q_2$}
\begin{tabular}{|l|c|c|c|c|c|c|}\hline
 & BDF(1) & BDF(2) & BDF(3) & BDF(4) & BDF(5) & BDF(6) \\ \hline
 DTOCF     & 0.4799 & 0.4061 & 0.3195 & 0.3097 & 0.2891 & 0.1871 \\ \hline
 APOCF     & 0.1000 & 0.0631 & 0.0631 & 0.0316 & 0.0316 & 0.0631 \\ \hline
 DPOCF(6)  & 0.5998 & 0.6641 & 0.6984 & 0.7475 & 0.7926 & 0.7880 \\ \hline
 DPOCF(12) & 0.5231 & 0.5815 & 0.5894 & 0.6165 & 0.6703 & 0.6516 \\ \hline
 DPOCF(20) & 0.5231 & 0.5149 & 0.4862 & 0.5044 & 0.5632 & 0.5326 \\ \hline
 DPOCF(30) & 0.4282 & 0.4144 & 0.3569 & 0.3636 & 0.3537 & 0.3136 \\ \hline
 DPOCF(40) & 0.3989 & 0.3649 & 0.2982 & 0.3021 & 0.2888 & 0.2398 \\ \hline
 DPOCF(60) & 0.3506 & 0.3021 & 0.2267 & 0.2267 & 0.2107 & 0.1672 \\ \hline
\end{tabular}}
\end{table}

\subsection{Accelerating Effect by Windowing Technique}
\label{sec-num-aewt}

In fact, the accelerating effect by windowing technique is revealed
in the sense of practical optimal convergence factor in subsection \ref{sec-num-ocf}.
The larger the number
of windows, the smaller the practical optimal convergence factor, or equivalently, the
faster the convergence rate of the DABSOR method. In this subsection, the accelerating
effect by windowing technique is exhibited in the sense of average number of
iterative steps of the DABSOR method on each window.

The comparisons of average number of iterative steps of the DABSOR method
based on six different linear multistep formulae like BDF(1-6),
two grid sizes as $12\times 12$, $24\times 24$ and two different choices of
preconditioners $Q$ in Table \ref{tab:caseQ} are outlined in Tables
\ref{tab:dwr_size12_Q1}-\ref{tab:dwr_size24_Q2}. In these tables, NoW stands for the
number of windows, NoU is the number of unknowns on each window, and ``---''
means the DABSOR method does not converge on at least one of the windows
in 800 iterative steps. Obviously, the average number of iterative steps
decreases in all kinds of situations when NoW
increases, which implies a decrease to the computation load.
In another word, the computation efficiency of the DABSOR method is
greatly improved by applying windowing technique. Moreover, the average
number of iterative steps of the DABSOR method based on preconditioner
$Q_2$ is apparently smaller than that based on preconditioner $Q_1$, which tells
the same story as in subsections \ref{sec-num-occ} and \ref{sec-num-ocf}.

It is known that high order time stepping schemes lead to high accuracy
for solving ODEs and DAEs, meanwhile the computation load increases
intensively. It is found in Tables \ref{tab:dwr_size12_Q1}-\ref{tab:dwr_size24_Q2}
that the average number of iterative steps of the DABSOR method increases
when the order of BDF becomes higher. For extreme situations, when NoW
is small, the DABSOR methods based on high order BDF
methods do not converge in 800 iterative steps. Hence, there should
be a balance between computation accuracy and computation efficiency.

\begin{table}
\setlength{\abovecaptionskip}{0pt}
\setlength{\belowcaptionskip}{10pt}
\centering{
\caption{\label{tab:dwr_size12_Q1} Average number of iteration steps of the DABSOR method
based on $12\times 12$ grid and preconditioner $Q_1$}
\begin{tabular}{|c|c|c|c|c|c|c|c|}\hline
NoW & BDF(1) & BDF(2) & BDF(3) & BDF(4) & BDF(5) & BDF(6) & NoU\\ \hline
 1 & 161.0 & 412.0 & ---   & ---   & ---   & ---    & 51840\\ \hline
 2 & 91.5  & 171.0 & 392.0 & 677.5 & ---   & ---    & 25920\\ \hline
 3 & 68.0  & 112.7 & 233.0 & 369.0 & 472.7 & 691.7  & 17280\\ \hline
 4 & 58.8  & 77.2  & 160.0 & 200.5 & 294.2 & 314.0  & 12960\\ \hline
 5 & 48.8  & 64.8  & 118.8 & 147.0 & 201.0 & 211.8  & 10368\\ \hline
 6 & 42.8  & 58.0  & 92.0  & 126.5 & 146.2 & 152.5  & 8640\\ \hline
\end{tabular}}
\end{table}
\begin{table}
\setlength{\abovecaptionskip}{0pt}
\setlength{\belowcaptionskip}{10pt}
\centering{
\caption{\label{tab:dwr_size12_Q2} Average number of iteration steps of the DABSOR method
based on $12\times 12$ grid and preconditioner $Q_2$}
\begin{tabular}{|c|c|c|c|c|c|c|c|}\hline
NoW & BDF(1) & BDF(2) & BDF(3) & BDF(4) & BDF(5) & BDF(6) & NoU\\ \hline
 1 & 160.0 & 428.0 & 697.0 & ---   & ---   & ---   & 51840\\ \hline
 2 & 76.0  & 134.0 & 210.5 & 373.5 & 717.5 & ---   & 25920\\ \hline
 3 & 49.7  & 75.0  & 107.7 & 159.3 & 193.3 & 264.3 & 17280\\ \hline
 4 & 36.8  & 47.5  & 71.2  & 87.5  & 121.2 & 125.2 & 12960\\ \hline
 5 & 31.2  & 35.0  & 43.6  & 59.8  & 67.6  & 83.2  & 10368\\ \hline
 6 & 27.7  & 30.7  & 35.3  & 34.7  & 47.7  & 47.7  & 8640\\ \hline
\end{tabular}}
\end{table}
\begin{table}
\setlength{\abovecaptionskip}{0pt}
\setlength{\belowcaptionskip}{10pt}
\centering{
\caption{\label{tab:dwr_size24_Q1} Average number of iteration steps of the DABSOR method
based on $24\times 24$ grid and preconditioner $Q_1$}
\begin{tabular}{|c|c|c|c|c|c|c|c|}\hline
NoW & BDF(1) & BDF(2) & BDF(3) & BDF(4) & BDF(5) & BDF(6) & NoU\\ \hline
 1 & 154.0 & 350.0 & 721.0 & ---   & ---   & ---   & 207360 \\ \hline
 2 & 75.0  & 129.0 & 261.0 & 363.0 & 719.5 & ---   & 103680 \\ \hline
 3 & 54.7  & 80.7  & 162.0 & 176.3 & 283.3 & 379.0 & 69120  \\ \hline
 4 & 46.5  & 62.8  & 115.2 & 122.2 & 209.8 & 189.0 & 51840  \\ \hline
 5 & 39.0  & 53.8  & 87.8  & 91.6  & 148.4 & 129.6 & 41472  \\ \hline
 6 & 34.2  & 47.5  & 69.2  & 69.2  & 110.8 & 93.5  & 34560  \\ \hline
\end{tabular}}
\end{table}
\begin{table}
\setlength{\abovecaptionskip}{0pt}
\setlength{\belowcaptionskip}{10pt}
\centering{
\caption{\label{tab:dwr_size24_Q2} Average number of iteration steps of the DABSOR method
based on $24\times 24$ grid and preconditioner $Q_2$}
\begin{tabular}{|c|c|c|c|c|c|c|c|}\hline
NoW & BDF(1) & BDF(2) & BDF(3) & BDF(4) & BDF(5) & BDF(6) & NoU\\ \hline
 1 & 168.0 & 366.0 & ---   & ---   & ---   & ---   & 207360 \\ \hline
 2 & 76.0  & 138.5 & 238.0 & 325.5 & 601.0 & ---   & 103680 \\ \hline
 3 & 48.0  & 75.7  & 110.7 & 150.0 & 226.0 & 338.0 & 69120  \\ \hline
 4 & 35.0  & 47.5  & 68.2  & 98.0  & 118.8 & 153.8 & 51840  \\ \hline
 5 & 29.6  & 36.8  & 49.0  & 62.4  & 83.2  & 102.4 & 41472  \\ \hline
 6 & 25.7  & 33.0  & 38.2  & 47.5  & 60.5  & 60.7  & 34560  \\ \hline
\end{tabular}}
\end{table}

\section{Concluding Remarks}
\label{sec-cremarks}
This paper studies the general theory of the discrete-time waveform relaxation methods.
Then, the DABSOR method is proposed for solving linear DAEs (\ref{dis-dae-LCC})
derived from time-dependent Stokes equations (\ref{time-stokes}). The convergence property
and optimality of the DABSOR method are stated in detail. Due to the requirement of
acceleration, the DABSOR method is integrated with windowing technique, which leads to
great acceleration as shown in Section \ref{sec-num}. In fact, further acceleration by
algebraic techniques like Krylov subspace on each
window is another path to improve the DABSOR method. The future work will follow right this path.

\end{document}